\documentclass[11pt]{amsart}

\usepackage{amssymb} 
\usepackage[latin1]{inputenc}
\usepackage[matrix,arrow,curve]{xy}
\usepackage[mathscr]{eucal}
\usepackage{epic,eepic}
\usepackage{bbm}

\makeindex
\newtheorem{itheorem}{Theorem}
\newtheorem{iprop}{Proposition}
\newtheorem{theorem}{Theorem}[subsection]
\newtheorem{lemma}[theorem]{Lemma}
\newtheorem{scholium}[theorem]{Scholium}
\newtheorem{prop}[theorem]{Proposition}

\newtheorem{cor}[theorem]{Corollary}
\theoremstyle{definition}
\newtheorem{definition}[theorem]{Definition}
\newtheorem{remark}[theorem]{Remark}
\newtheorem{example}[theorem]{Example}

\newtheorem{question}[theorem]{Question}

\newtheoremstyle{numero}
{1ex}
{1ex}
{\it}
{1cm}
{\scshape}
{}
{4pt}
{}
\theoremstyle{numero}

\newenvironment{pf}
{\medskip\noindent {\it Proof --- \ }}
{\hfill\nobreak $\Box$ \par\bigbreak}

\renewcommand{\div}{{\text{div}}}

\renewcommand{\P}{{\mathbb P}}

\newcommand{\F}{ \mathbb F} 
\newcommand{\C}{{ \mathbb C  }}
\newcommand{\R}{{ \mathbb R  }}
\newcommand{\Q}{{ \mathbb Q } }
\newcommand{\Z}{{ \mathbb Z  }}

\renewcommand{\ker}{{\text{Ker}\,}}
\newcommand{\Gal}{{\mathrm{Gal}\,}}

\newcommand{\anneau}{{ \mathcal O}}

\newcommand{\Id}{{\text{Id}}}

\newcommand{\PSL}{{\text{PSL}}}

\newcommand{\spec}{{\text{Spec\,}}}

\newcommand{\Pic}{{\text{\rm Pic}}}

\newcommand{\pr}{{\text{pr}}}

\newcommand{\Fc}{{\mathcal{F}}}
\newcommand{\Fcb}{{\mathcal{F}b}}

\newcommand{\ch}{{\text{char\,}}}

\newcommand{\ord}{{\text{ord}}}
\newcommand{\tr}{{\text{tr\,}}}

\renewcommand{\hat}{\widehat}

\renewcommand{\Gal}{{\rm Gal}}

\newcommand{\CC}{{\mathcal{C}}}
\newcommand{\mi}{{\text{min}}}
\newcommand{\sep}{{\text{sep}}}
\newcommand{\Div}{{\text{Div}}}
\newcommand{\perf}{{\text{perf}}}
\newcommand{\np}{{\rm{np}}}
\newcommand{\crit}{{\rm{crit}}}

\newcommand{\Eb}{E_{\text{backward}}}

\newcommand{\Kb}{K_{\text{backward}}}
\newcommand{\poincare}{{\mathcal{H}}}

\newcommand{\Lc}{{\mathcal{L}}}
\begin{document}
\setcounter{tocdepth}{1}
\baselineskip 15pt
\title{On self-correspondences on curves}
\author{Joël Bellaïche}
\address{MS 050 \\ Brandeis University \\ 45 South Street \\ Waltham, MA \\ 02453 \\USA}
\email{jbellaic@brandeis.edu}
\begin{abstract}
We study the algebraic dynamics of self-correspondences on a curve. A self-correspondence on a (proper and smooth) curve $C$ over an algebraically closed field is the data of another curve $D$ and two non-constant separable morphisms $\pi_1$ and $\pi_2$ from $D$ to $C$. A subset $S$ of $C$ is {\it complete} if $\pi_1^{-1}(S)=\pi_2^{-1}(S)$. We show that self-correspondences are divided into two classes: those that have only finitely many finite complete sets, and those for which $C$ is a union of finite complete sets. The latter ones are called {\it finitary}, happen only when $\deg \pi_1=\deg\pi_2$ and have a trivial dynamics. For a non-finitary self-correspondence in characteristic zero, we give a sharp bound for the number of étale finite complete sets.
\end{abstract}
\maketitle

\section*{Introduction}

Let $k$ be a field, and let  $C$ a be smooth, proper and geometrically irreducible curve over $k$. 
By a {\it self-correspondence}\footnote{We adopt the definition of \cite{BP} and \cite{raju}. In part of the literature, a self-correspondence is defined instead as a divisor in the surface $C \times C$. The two notions are equivalent. To get a divisor of $C\times C$ using our definition, take $(\pi_1 \times \pi_2)(D)$, with multiplicities if several components of $D$ have the same image in $C \times C$. To get from a divisor $\Delta$ of $C \times C$ a self-correspondence according to our definition, 
take the union of the normalization of each component of $\Delta$ repeated according to multiplicity. Our definition makes clearer
the concepts of {\it étale} or {\it equiramified} complete sets, which is central in our study.} on $C$ (defined over $k$), we mean the data of a smooth and proper scheme $D$ over $k$, such that every connected component of $D$ is a geometrically irreducible curve, and two $k$-morphisms $\pi_1$ and $\pi_2$ from $D$ to $C$, non-constant and separable on every connected component of $D$. We denote by $(D,\pi_1,\pi_2)$, or often simply by $D$, that self-correspondence. 

Fixing an algebraic closure $\bar k$ of $k$, the intuitive way to think of a self-correspondence is as a multi-valued map from $C(\bar k)$ to itself defined by polynomial equations with coefficients in $k$, namely the map $x \mapsto \pi_2(\pi_1^{-1}(x))$. We call this multi-valued map the {\it forward map} of the self-correspondence $D$.

Self-correspondences generalize endomorphisms: given an endomorphism $f$ of a curve $C$, one can think of it as the self-correspondence $D_f := (C,\Id_C,f)$. Better than endomorphisms, a self-correspondence $D=(D,\pi_1,\pi_2)$ always has a  {\it transpose}, denoted by ${}^t D$ and defined by ${}^t D = (D,\pi_2,\pi_1)$. The forward-map of ${}^t D$, namely $x \mapsto \pi_1(\pi_2^{-1}(x))$ is called the {\it backward map} of $D$. 

Let us introduce our fundamental terminology. A {\it forward-complete} (resp. {\it backward-complete}) set is a subset $S$ of $C(\bar k)$ that is stable by the forward (resp. backward) map (i.e. $\pi_1^{-1}(S) \subset \pi_2^{-1}(S)$, resp. $\pi_2^{-1}(S) \subset \pi_1^{-1}(S)$), and a {\it complete} set if a set which is both backward and forward-complete. An {\it irreducible complete set} is a minimal non-empty complete set. For $z \in D(\bar k)$ we write $e_{1,z}$, $e_{2,z}$ for the ramification index of $\pi_1$ and $\pi_2$ at $z$ and we say that $z \in D(\bar k)$ is {\it ramification increasing} if  $e_{1,z} \leq e_{2,z}$, {\it ramification decreasing} if $e_{1,z} \geq e_{2,z}$, {\it equiramified} if $e_{1,z}=e_{2,z}$ and {\it étale} if $e_{1,z}=e_{2,z}=1$. All points of $D(\bar k)$ except possibly a finite number are étale, but important phenomena occur at non-étale point as well. We say that a subset $S \subset C(\bar k)$ is  {\it ramification increasing} (resp. {\it ramification decreasing, equiramified, étale}) if every $z \in D(\bar k)$ such that $\pi_1(z) \in S$ and $\pi_2(z) \in S$ is ramification increasing (resp. etc.).
\par \bigskip
The aim of this article is to answer elementary questions about finite complete sets for self-correspondences, such as when can there be infinitely many finite complete sets? This is a basic and fundamental question on the dynamics of self-correspondence. There is a relatively extensive literature on the subject of dynamics of self-correspondences on curves and even, lest often, algebraic varieties. Most of its literature is concerned about correspondences over the complex numbers, see for instance \cite{fatou}, \cite{B0}, \cite{B1}, \cite{BP0}, \cite{BP}, \cite{BL}, \cite{BS} (on $\P^1$), \cite{Dinh} (on $\P^k$), \cite{DS} (on general varieties), \cite{DKW} (on general curves), but also over number fields, see \cite{A}, \cite{I1}, \cite{I2}, finite fields
 \cite{HP} and general fields \cite{TTT}). To the best of our knowledge, the question we have in mind has not been solved, and not even asked. A partial exception  is the recent article of Raju (\cite{raju}), essentially the second chapter of his PhD thesis at Columbia University. Though the focus of the paper is different, as Raju is concerned with general correspondences rather than self-correspondences, we borrow several ideas and concepts to him, and we gladly acknowledge our debt to his work.

\par \medskip
Our first main result concerning finite complete sets is the following (see Theorem~\ref{finitary}) :
\begin{itheorem} A self-correspondence $(D,\pi_1,\pi_2)$ on a curve $C$ over $k$ has infinitely many finite complete sets if and only if there exists a non-constant $k$-morphism $f:C \rightarrow \P^1_k$ such that $f \circ \pi_1= f \circ \pi_2$.
\end{itheorem}
The method of proof of Theorem 1 is number-theoretic. Specifically, we use the  famous theorem of Mordel-Weil-Néron asserting that the group of rational points of an abelian variety over a finitely generated field is a finitely generated abelian group. Theorem 1 seems difficult to prove by purely algebraico-geometric methods or by complex-analytic ones in the case $k=\C$. 

 As a trivial consequence of Theorem 1, one sees that as soon as a self-correspondence has infinitely many finite complete sets, then in fact all its irreducible complete sets are finite, and moreover they have cardinality bounded by some integer $M$. A self-correspondence satisfying this property will be said {\it finitary}; its dynamical study is essentially trivial, in the sense that it reduces to dynamical questions over finite sets. 
 
It is clear  that only self-correspondences for which $\deg \pi_1=\deg \pi_2$ (such a self-correspondence is said to be {\it balanced}) can be finitary.
This explains why the notion of being finitary does not appear in the classical theory of complex dynamics, where one consider endomorphisms of $\P^1$, which as correspondences have $(\deg \pi_1,\deg \pi_2)=(1,d)$, the case $d=1$ being excluded as trivial.  Even 
among balanced self-correspondences, the notion of finitary self-correspondence is very restrictive. They are in practice exceptions to all interesting general statements about the
dynamics of self-correspondences. Non-finitary self-correspondences are the natural domain of study of the dynamics of self-correspondences. For instance,
an interesting result on the existence of a canonical invariant measure for a balanced self-correspondence on a curve over $\C$ has recently been proved by 
Dinh, Kaufmann and Wu (\cite{DKW}) but only under a quite restrictive condition on the self-correspondence (see Remark~\ref{remarkDKW}). We believe that their main result (that the iterated pull-back of every smooth measure converges to the canonical measure) holds for all non-finitary correspondences, and we plan to come back to this question in a subsequent work.

\par \medskip

For a non-balanced correspondence, there are no étale finite complete set, and finitely many non-étale ones. Moreover one can give an upper bound (in terms of the genera and degrees of the curves and maps involved)
on the size of their union (see Prop.~\ref{unbalancedprop}). 

Balanced non-finitary correspondence are much more subtle: they may have (finitely many) étale and non-étale finite complete sets; it is easy to give a bound on the number of  finite non-étale complete sets, but we do not know how to bound their size. Our main objective is to bound the number of finite étale complete sets. With methods similar to those of the proof of Theorem 1, we are able to offer a bound (which happens to be optimal) only in some specific cases: when $k$ is algebraic over a finite field (see Prop.~\ref{propetalefinite}); when $k$ is arbitrary but $C=\P^1_k$ (see Prop.~\ref{caseP1}); and two other results when $D$ is symmetric, that is $D \simeq {}^t D$ (see Propositions~\ref{propraju} and~\ref{mochi}). 
But for more general results we need different, operator-theoretic. methods.
\par \medskip
A self-correspondence $D$ over $C$ defines in an natural way a $k$-linear endomorphism $T_D$ of the field of rational functions $k(C)$ of $C$, see \S\ref{defTD}.
Whenever $S$ is a forward-complete set, $T_D$ stabilizes the subring $B_S$ of $k(C)$ of functions whose all poles are in $S$.
If moreover $S$ is ramification-increasing, $T_D$ stabilizes as well the natural filtration $(B_{S,n})_{n \geq 0}$, ``by the order of the poles"
of $B_S$. 

The dynamical study of that action of $T_D$ on the filtered ring $B_S$, when $S$ is in particular étale complete, was the original motivation of this work. In fact, Hecke operators appearing in the theory of modular forms are of this type. Surprising results concerning the dynamics of the operators $T_D$
for self-correspondences over finite fields have been obtained in a recent work by Medvedovsky (\cite{medvedovsky}), and, applied to Hecke operators, those results provide a new and elementary proof of certain deep modularity result of Gouvêa-Mazur (\cite{GM}). We plan to come back to these questions on a subsequent work. But in this paper, we content ourselves  to use the operators $T_D$ to obtain new informations on the dynamics of $D$,

We say that $D$ is {\it linearly finitary} if there is a monic polynomial $Q$ in $k[X]$ such that $Q(T_D)=0$ as endomorphisms of $k(C)$. 
That $D$ is linearly finitary means that the dynamics of $T_D$ is trivial, in the sense that it is similar to the one of an operator on a finite-dimensional vector space. We prove the following:

\begin{iprop} A finitary self-correspondence is  linearly finitary. The converse is true in characteristic zero.  \end{iprop}
The direct sense is very easy. We prove the converse using graph-theoretic methods and Theorem 1 (see Prop.~\ref{linfinfin}.)

Our second main result is the following (see Theorem~\ref{notthree} which is slightly more precise) :
\begin{itheorem} If a self-correspondence $D$ has three irreducible étale finite complete sets, then it is linearly finitary. In particular, in characteristic zero, 
a non-finitary correspondence has at most two irreducible finite étale complete sets.
\end{itheorem}
We give a brief description of the idea of the proof, which uses only elementary methods : the theorem of Riemann-Roch and linear algebra. The first étale finite  complete set $S$ is used to define, as above, the natural filtration $(B_{S,n})_{n \geq 0}$, ``by the order of the poles" of $B_S$, which is stabilized by $T_D$. Is this filtration split as a $k[T_D]$-filtration? In general, this has no reason to be true. But with a second irreducible étale complete set $S'$, one can show that this filtration is ``almost split'', namely that there exists for every $n$ a $T_D$-stable subspace $V_{S,S',n}$ in $B_{S,n+1}$ such that $B_{S,n} + V_{S,S',n} = B_{S,n+1}$ and $\dim V_{S,S',n}$ bounded\footnote{The filtration $(B_{S,n})$ of $B_S$ would be $T_D$-split if in addition of $B_{S,n} + V_{S,S',n} = B_{S,n+1}$ we required $B_{S,n} \cap V_{S,S',n} = 0$, or equivalently
$\dim V_{S,S',n} = \dim B_{S,n+1}/ B_{S,n}$.  One has  $\dim B_{S,n+1}/B_{S,n} \leq |S|$ for every $n$ (with equality for $n$ large enough), so requiring that $\dim V_{S,S',n}$ is bounded independently of $n$ is a qualitative version of that property. Hence the phrase ``almost split".} independently of $n$. We prove this by defining $V_{S,S',n}$ as the space of functions in $B_{S,n+1}$ that vanishes on $S'$ at a suitable order, and using Riemann-Roch.
  Now a third finite étale complete set $S''$ give a second quasi-splitting $V_{S,S'',n}$ of the filtration $B_{S,n}$. Using Riemann-Roch, we prove that the two filtrations are orthogonal, in the sense that $V_{S,S',n} \cap V_{S,S'',n}=0$ for $n$ large enough. Then, a linear algebra argument shows that all eigenvalues of $T_D$ appearing in $B_{S,n+1}$ (for $n$ large enough) already appear in $B_{S,n}$ and Theorem 2 follows.

\tableofcontents

\section{Self-correspondences}

\subsection{Curves}
Let $k$ be a field. By a {\it curve} over $k$ we shall mean a non-empty proper and smooth scheme over $\spec k$ which is equidimensional of dimension $1$ and geometrically connected.

If $C$ is a curve over $k$, we denote by $k(C)$ the function field of $C$. If $C$ and $C'$ are two curves over $k$, a non-constant morphism of $k$-schemes $\pi: C \rightarrow C'$ is finite and flat, 
hence surjective  and thus defines a morphism of $k$-extensions $\pi^\ast: k(C') \rightarrow k(C)$, which makes $k(C)$ a finite extension of $k(C')$. Explicitly, if $f \in k(C)$ is seen as a morphism from $C$ to $\P^1$, $\pi^\ast f = f \circ \pi$. We say that $\pi$ is {\it  separable} if $k(C)$ is a separable extension of $k(C')$.


\subsection{Self-correspondences}

Given a curve $C$ over $k$, a {\it self-correspondence} $(D,\pi_1,\pi_2)$ on  $C$ is the data of a noetherian reduced $k$-scheme $D$ whose connected components are curves $D_i$ over $k$, with two morphisms $\pi_1$ and $\pi_2$ from $D$ to $C$ whose restrictions to each $D_i$ are non-constant and separable.
Often, the morphisms $\pi_1$ and $\pi_2$ will be implicit and we shall simply denote by $D$ the self-correspondence $(D,\pi_1,\pi_2)$


Self-correspondences on a fixed curve $C$ over $k$ naturally form a category: a morphism from $(D,\pi_1,\pi_2)$ to $(D',\pi'_1,\pi'_2)$ is a surjective $k$-morphism $h: D \rightarrow D'$ such that $\pi'_i \circ h = \pi_i$ for $i=1,2$. In particular we have a notion of {\it isomorphism} of self-correspondences.

\begin{example} \label{hecke} Let $C$ be the complete Igusa curves of level $N$ (with $(N,p)=1$) over $\F_p$. Recall (\cite[pages 460-462]{gross}, where the curve is denoted by $I_1(N)$) that $C$ is the smooth completion of the affine Igusa curve, which is defined as the  moduli space for triples $(E,\alpha,\beta)$, where $E$ is an elliptic curve over a scheme of characteristic $p$, $\alpha$ an embedding $\mu_N \hookrightarrow E$ and $\beta$ an embedding $\mu_p \hookrightarrow E$.

For $l$ a prime number not dividing $Np$, we define the {\it Hecke correspondence} $D_l$, moduli space for quadruples $(E,\alpha,\beta,H)$ where $(E,\alpha,\beta)$ is as above and $H$ is a subgroup scheme of $E$  locally  of order $l$. Define $\pi_1: D_l \rightarrow C$ as just forgetting $H$, and $\pi_2$ as sending $(E,\alpha,\beta,H)$ to $(E/H,\alpha',\beta')$ where $\alpha'$ and $\beta'$ are defined in the obvious manner. Then $(D_l,\pi_1,\pi_2)$ is a self-correspondence, called the Hecke correspondence at $l$, on the Igusa curve $C$.
\end{example}

\subsection{Bi-degree}
 
If $D$ is a finite disjoint union of curves $D=\coprod_i D_i$, $C$ a curve, and $\pi:D \rightarrow C$ a map non-constant on every component $D_i$ of $D$, then we define $\deg \pi$ as $\sum_{i} \deg \pi_{|D_i}$.

The {\it bi-degree} of a self-correspondence $(D,\pi_1,\pi_2)$ on a curve $C$ is the ordered pair of integers $(\deg \pi_1,\deg \pi_2)$. It is often denoted $(d_1,d_2)$. A self-correspondence is {\it balanced} when $d_1=d_2$. For example, the bi-degree of the Hecke correspondence $D_l$ of Example~\ref{hecke} is $(l+1,l+1)$. 

\subsection{Transpose}
  
If $(D,\pi_1,\pi_2)$ is a self-correspondence on $C$, so is $(D,\pi_2,\pi_1)$, called the {\it transpose}  of $(D,\pi_1,\pi_2)$.
We shall denote this correspondence by ${}^t D$. Its bi-degree is $(d_2,d_1)$, if the degree of $D$ is $(d_1,d_2)$.

\subsection{Self-correspondence of morphism type} Let $f$ be a non-constant separable morphism from $C$ to $C$. We denote by  $D_f$ the self-correspondence $(C,\Id_C,f)$. A self-correspondence $D$ over $C$ is {\it of morphism type} if it is isomorphic to some $D_f$. Equivalently, $D$ is of morphism type if and only if its bi-degree has the form $(1,d)$. Thus we see that the transpose of a self-correspondence of the form $D_f$ is of morphism type only when $f$ is an isomorphism, and in this case ${}^t D_f \simeq D_{f^{-1}}$.

\subsection{Minimal self-correspondences}  A self-correspondence is {\it minimal} if the map $\pi_1 \times \pi_2: D \rightarrow C^2$
is generically injective, that is there are only finitely many points of $C^2$ such that the fiber of $\pi_1 \times \pi_2$ has more than one
element. Equivalently, $D$ is minimal if $k(D)$ is generated, as a $k$-algebra, by its two subfields $\pi_1^\ast(k(C))$ and $\pi_2^\ast(k(C))$.

For any self-correspondence $D$ on $C$, there is a unique pair $(D_\mi,h)$ where $D_\mi$ is a minimal self-correspondence over $C$  and $h: D_\mi\rightarrow D$ a morphism of self-correspondences on $C$: take $D_\mi$ the normalization of the image of $D$ by $\pi_1 \times \pi_2$. 

Note that a minimal self-correspondence is rigid, i.e. has a trivial group of automorphism. 

\subsection{Symmetric self-correspondences}
 
 A self-correspondence $D$ is {\it symmetric} if ${}^t D \simeq D$. Obviously a symmetric self-correspondence is balanced.
 A self-correspondence is symmetric if and only if there exists an automorphism $\eta$ of the $k$-scheme $D$ such that $\pi_1 \circ \eta = \pi_2$. When $D$ is minimal, this automorphism $\eta$ is necessarily an involution, since $\eta^2$ is an automorphism of the self-correspondence $D$, which is rigid.

The Hecke correspondence $D_l$ on the Igusa curve (see Example~\ref{hecke}) is symmetric with $\eta(E,\alpha,\beta,H)=(E/H,\alpha',\beta',E[l]/H)$.

\subsection{Terminology on directed graphs}

\label{terminologydirectedgraph}

By a {\it directed graph} we shall mean the data $\Gamma = (V,Z,s,t)$ of two sets $V$ and $Z$ and two maps $s,t: Z \rightarrow V$. Elements of $V$ 
are called {\it vertices}, elements of $Z$ are called {\it edges}, and for $z \in Z$, $s(z)$ is the {\it source} and $t(z)$ the {\it target} of $z$.
In particular, {\it self-loops} (i.e. edges $z$ such that $s(z)=t(z)$) and {\it repeated edges} (i.e. edges $z_1 \neq z_2$ such that $s(z_1)=s(z_2)$ and $t(z_1)=t(z_2)$) are allowed.

We use the usual terminology for a directed graph: a {\it forward-neighbor} (resp. {\it backward-neighbor}) of a vertex $x$ is a vertex $y$ such that there is an edge $z$ with source $x$ and target $y$ (resp. with source $y$ and target $x$). More generally, we say 
that $y$ is a {\it $k$-forward-neighbor} of $x$ if $y=x$ when $k=0$, or if $y$ is a forward-neighbor of a $k-1$-forward-neighbor when $k \geq 1$. 

A subset $S$ of $V$ is said to be {\it forward-complete} (resp. {\it backward-complete}) if it contains every forward-neighbors (resp. backward-neighbors) of its vertices. A subset $S$ of $V$ is {\it complete} if it is both backward-complete and forward-complete. A complete subset $S$ is {\it irreducible} if it is non-empty and has no complete proper non-empty subset. 

It is clear that the irreducible complete subsets of a directed graph are the connected components of the non-directed graph it defines, and that the complete sets are the union of connected components. A union and intersection of complete sets is complete, as is the complement of any complete set. Every complete set is a disjoint union of irreducible complete sets.

If $x,y$ are in $V$, a {\it directed path of length $n$} from $x$ to $y$ is a sequences of $n$ edges $p=(z_1,\dots,z_n)$ such that $s(z_1)=x$, $t(z_n)=y$ and $t(z_i)=s(z_{i+1})$ for $i=1,\dots,n-1$. A directed path from $x$ to $x$ is called a {\it directed cycle}.We shall denote by $\np_{x,y,n}$ the number of directed paths from $x$ to $y$.

If $k$ is a ring, and $\Gamma = (V,Z,s,t)$ is a directed graph, we define the {\it adjacency operator} of $\Gamma$, $A_\Gamma : \CC(V,k) \rightarrow \CC(V,k)$ on the $k$-module $\CC(V,k)$ of maps from $V$ to $k$, by the formula 
$$(A_\Gamma f)(y) = \sum_{z \in Z, t(z)=y}  f(s(z)).$$ The matrix of $A_\Gamma$ in the canonical basis of $\CC(V,k)$ is the {\it adjacency matrix} of $\Gamma$.
By induction, we check that if $x,y\in V$ and $n \geq 1$ an integer, then 
\begin{eqnarray} \label{Agn}  (A_\Gamma^n(\delta_x)) (y) = \np_{x,y,n}.\end{eqnarray}

We shall sometimes consider functions $f: V \rightarrow k \cup \{\infty\}$ where $\infty$ is a symbol not in $k$. For such a function, we define $A_\Gamma f$ by the same formula as above, with the convention that the sum of the right hand-side is $\infty$ if exactly one of its term is $\infty$, and is undefined if two or more of its terms are $\infty$. 


\subsection{The directed graph attached to a self-correspondence}

If $(D,\pi_1,\pi_2)$ is a self-correspondence of $C$ over $k$, and $\bar k$ is a fixed algebraic closure of $k$,
 we define the oriented graph $\Gamma_D$ attached to $D$ as $(C(\bar k),D(\bar k),\pi_1,\pi_2)$.
 
 If $z \in D(\bar k)$ is an edge, we write  $e_{i,z}$ for the index of ramification of $\pi_i$ at $z$.
  An edge $z \in D(\bar k)$ is said {\it ramification increasing} (resp. {\it equiramified}, resp. {\it étale}) if $e_{1,z} \leq e_{2,z}$ (resp. $e_{1,z}=e_{2,z}$, resp. $e_{1.z}=e_{2,z}=1$). We observe that there only finitely many edges that are not étale (and a fortiori, not equiramified or ramification increasing). This is because $\pi_1$ and $\pi_2$ are assumed separable. 
 
A subset $S$ of $C(\bar k)$ is said {\it ramification increasing}, {\it equiramified}, {\it étale} if all the edges whose both source and target are in $S$ are {\it ramification increasing}, etc.
 
For any vertex $x \in C(\bar k)$, we have the formula
\begin{eqnarray} \label{sume1} \sum_{z \in D(\bar k), \pi_1(z)=x} e_{1,z} &=& d_1\\
\label{sume2} \sum_{z \in D(\bar k), \pi_2(z)=x} e_{2,z} &=& d_2. \end{eqnarray}
In particular, there are at most $d_1$ edges with source $x$ and $d_2$ edges with target $x$, and the directed graph $\Gamma_D$ is locally finite.
Given a finite set of vertices $S \subset C(\bar k)$, one has by summing the above formula:
\begin{eqnarray} \label{sumS1} \sum_{z \in \pi_1^{-1}(S)} e_{1,z} &=& d_1 |S|\\
\label{sumS2} \sum_{z \in \pi_2^{-1}(S)} e_{2,z} &=& d_2 |S|. \end{eqnarray}

\begin{remark} The directed graph of ${}^t D$ is the directed graph of $D$ with source and target maps  exchanged.\end{remark}

\begin{lemma} Let $k$ be a finitely generated extension of a prime field. Then $k$ has extensions of arbitrary large prime degrees.
\end{lemma}
\begin{pf} Let $\F$ be the prime subfield of $k$, and let $T_1,\dots,T_n$ be a transcendence basis of $k$ over $\F$.
Thus if $k_0=\F(T_1,\dots,T_n)$, $k$ has finite degree $d$ over $k_0$. The field $k_0$ admits extension of any prime degree $p$ : if $n=0$, then $k_0=\F_p$ or $\Q$ and the result is well-known, and if $n \geq 1$, for $p$ prime, the polynomial $X^p-T_1$ has no root in $k_0$ hence is irreducible over $k_0$ (see \cite[Theorem 9.1]{Lang}). If $p \nmid d$, the composition of an extension of degree $p$ of $k_0$ with $k$ is an extension of $k$ of degree $p$.
\end{pf}

\begin{prop} \label{infinitecomplete} The directed graph $\Gamma_D$ has infinitely many irreducible complete sets. All but finitely many of them are étale.
\end{prop}
\begin{pf}
Let us consider $C$ and $D$ as embedded in a projective space over $k$ (say $\P^3_k$),  and let $k_0$ be the subfield of $k$ generated over the prime subfield of by the coefficients of the projective equations of $C$ and $D$ and the coefficients of the polynomials defining $\pi_1$ and $\pi_2$. Replacing $k$ by $k_0$ we may assume that $k$ is of finite type over its prime subfield.

If $k$ is a finite type extension of its prime subfield, and $x \in C(k)$ is a vertex, then if $z \in D(\bar k)$ is an edge with source (resp. target) $x$, one has $z \in D(k')$ for $k'$ some finite extension of $k$ of degree $\leq d_1$ (resp. $\leq d_2$). Indeed, $z$ belongs to the schematic fiber of $\pi_1$ at $x$, which is a finite $k_1$-scheme of degree $d_1$. It follows that any forward-neighbor (resp. backward-neighbor) of $x \in C(k)$ is defined on an extension of degree $\leq d_1$ (resp. $\leq d_2$) of $k'$. By induction, any vertex in the same irreducible complete set as $x$ is defined over an extension of $k$ of degree whose all prime factors are $\leq \max(d_1,d_2)$.

Given a finite family $S_1,\dots,S_l$ of irreducible complete sets in $C(\bar k)$, pick points $x_1 \in S_1,\dots, x_l \in S_l$.
By replacing $k$ by a finite extension, we may assume that $x_1,\dots,x_l$ all belong to $C(k)$, and thus every points $x$ in $S_1 \cup \dots \cup S_l$ belong to $C(k')$ for $k'$ a finite extension of $k$ (depending on $x$) of some degree whose all prime factors are less than $\max(d_1,d_2)$.
By the above lemma, $k$ has  extensions of arbitrary large prime degrees, hence has an extension $k''$ of prime degree $p > \max(d_1,d_2)$ and $C$ has a point whose field of of definition contains $k''$. Such a point cannot belong to $S_1 \cup \dots \cup S_l$,
which shows that there are other irreducible complete sets in $C(\bar k)$. Therefore the number of irreducible complete sets is infinite.

 The second assertion is clear since there are only finitely many non-étale edges.
\end{pf}

\begin{example} \label{examplegraphHecke} The directed graph of the Hecke correspondence $D_l$ on the Igusa curve $C$ (see Example~\ref{hecke}) is well understood. Since $D_l$ is symmetric, we loose no information by forgetting the orientation of the edges and looking at $\Gamma_{D_l}$ as an undirected graph.

There are two obvious finite completes sets, the set of supersingular points and the sets of cusps. The complete set of supersingular points is étale (easy since $l \neq p$)  and irreducible (this can be proved by direct analysis, or, as Raju notes in \cite{raju}, simply as a consequence of Prop~\ref{propraju} below). The complete set of cusps may be reducible, and none of its irreducible components are étale (again a consequence of Prop.~\ref{propraju}). The other complete sets are all infinite and  have been called {\it isogeny volcanoes}: they consist in a cycle of some order $n$, with an infinite tree of valence $l+1$ attached to each vertices of that cycle. See \cite{sutherland}, \cite{kohel}.
\end{example}

\subsection{Sum and composition of self-correspondences}

 Let $(D,\pi_1,\pi_2)$ and $(D',\pi'_1.\pi'_2)$ be self-correspondences on a curve $C$, of bi-degree $(d_1,d_2)$ and $(d'_1,d'_2)$.
 
The {\it sum} of $D$ and $D'$, denoted by $D+D'$ is by definition the self-correspondence 
 $(D \coprod D',\pi_1 \coprod \pi'_1 ,\pi_2 \coprod \pi'_2)$ on $C$. It is obvious that the oriented graph $\Gamma_{D + D'}$ has the same vertices as $\Gamma_D$ and $\Gamma_{D'}$ and for set of edges the disjoint union of their set of edges. The bi-degree of $D + D'$ is $(d_1+d'_1,d_2+d'_2)$.
 
We define $D' \circ D$ as the scheme $D \times_{\pi_2,C,\pi'_1} D'$. The two projections $\pr_1$ and $\pr_2$ of this fibered product over $D$ and $D'$ are finite and flat, and its total ring of fractions is étale over $k(C)$, since it is the tensor product of the two separable extensions $k(D)$ and $k(D')$ over $k(C)$ (they are seen as extensions of $k(C)$ through $\pi_2^\ast$ and ${\pi'}_1^\ast$ respectively).
 In particular, $D' \circ D$ is proper over $k$, and all of its irreducible component have dimension $1$.
  We denote by $D'D$ the normalization of the reduced scheme attached to $D' \circ D$, and by $n :D'D \rightarrow D \circ D'$ the natural map. 
 
 $$
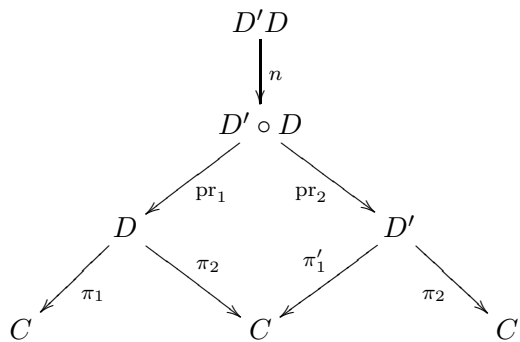
\xymatrix{ & & D'D \ar[d]^n & & \\ & & D' \circ D \ar[dl]^{\pr_1} \ar[dr]_{\pr_2} & & \\
 & D \ar[dr]^{\pi_2} \ar[dl]^{\pi_1} & & D' \ar[dl]_{\pi'_1} \ar[dr]_{\pi_2}& \\
 C & & C & & C }$$ 
 
Thus $D'D$ is a proper and smooth scheme of dimension $1$, that is a disjoint union of curves. 
Since $n$ is surjective, the restriction of $ \pr_1 \circ n$ and $ \pr_2 \circ n$ to every connected component of $D'D$ are surjective onto $D$ and $D'$ respectively and they are separable since $n$ induces an isomorphism on the total rings of fractions. Thus, $(D'D,\pi_1 \circ \pr_1 \circ n, \pi_2 \circ \pr_2 \circ n)$ is a self-correspondence on $C$, which we shall call {\it the composition of $D'$ with $D$}. Its bi-degree is $(d_1d'_1,d_2,d'_2)$.

\begin{example} If $f,g$ are morphisms from $C$ to $C$, then $D_f D_g = D_{fg}$. \end{example}

The directed graph $\Gamma_{D' \circ D}$ is  the graph whose vertices are those of $\Gamma_D$ or $\Gamma_{D'}$, and edges  $(z,z')$ where $z$ is an edge of $\Gamma_D$, $z'$ is an edge of $\Gamma_{D'}$ such that the source of $z'$ is the target of $z$. The source (resp. target) of the edge $(z,z')$ in $\Gamma_{D' \circ D}$ is the source of $z$ in $\Gamma_D$ (resp. the target of $z'$ in $\Gamma_{D'}$). Hence 
\begin{eqnarray} \label{AGammaDD} A_{\Gamma_{D' \circ D}} = A_{\Gamma_{D'}} \circ A_{\Gamma_D}. \end{eqnarray}

Since $n$ is surjective, generically an isomorphism, the directed graph $\Gamma_{D'D}$ is that of $\Gamma_{D' \circ D}$ described above, with finitely many new edges added, all those new edges having the same source and target that an already existing edge in $\Gamma_{D' \circ D}$.

\begin{lemma} \label{lemmaSDD'} If $S$ is a forward-complete (resp. backward-complete, resp. complete) set for $D$ and $D'$ then it is forward-complete (resp. etc.) for $D'D$. 

If $S$ is a complete étale set for $D$ and $D'$, then the restrictions of the  directed graphs of $D' \circ D$ and of $D'D$ to $S$ coincide, $S$ is also étale for $D'D$ and $A_{\Gamma_{D'D}}$ coincides with  $A_{\Gamma_{D'}} \circ A_{\Gamma_D}$ on $\CC(S,k)$.
 \end{lemma}
\begin{pf} The first assertion follows from what was said above. For the second, if $x \in S$, the fibers of $D$ and $D'$ at $x$ are both étale, and so is their tensor product $D_x \times_k D'_x$, which is $(D' \circ D)_x$, so $D' \circ D$ is étale at points above $x$. But then it is smooth, and $n: D'D \rightarrow D'\circ D$ is thus an isomorphism on some neighborhood of the points above $x$, so $D'D$ is also étale over $x$, and the last assertion follows from (\ref{AGammaDD}).
\end{pf}

It is clear that the composition of self-correspondences is associative up to obvious canonical identifications. 
We denote by $D^n$ the composition of $n$ copies of  the self-correspondence $D$.

\section{Finite complete sets}

\subsection{The unbalanced case}

\begin{prop}  \label{unbalancedprop} Let $(D,\pi_1,\pi_2)$ be a self-correspondence over $C$ of bi-degree $(d_1,d_2)$ with $d_1 < d_2$. Denote by $g_D$, $g_C$ the genera\footnote{The genus $g_D$ of a finite disjoint union of curves $D=\coprod_i D_i$ is defined here as the sum of the genera of $D_i$.
With the definition, and that of degree given in \S\ref{bideg}, Hurwitz formula is still valid for a map $\pi: D \rightarrow C$.}
of $D$ and $C$.
Then any finite backward-complete set $S$ satisfies $$|S| \leq 2 \frac{g_D - d_2 g_C +d_2 -1}{d_2-d_1}.$$ Moreover if such a set $S$ is ramification-increasing, then it is empty.
 \end{prop}
 \begin{pf} 
For the first assertion,
\begin{eqnarray*} |S| d_2 &=& \sum_{z \in \pi_2^{-1}(S)} e_{2,z} \text{ (by formula~(\ref{sumS2})})\\
&=& |\pi_2^{-1}(S)| +  \sum_{z \in \pi_2^{-1}(S)} (e_{2,z}-1)\\
&\leq& |\pi_2^{-1}(S)| +   \sum_{z \in D(\bar k)} (e_{2,z}-1) \\
&\leq& |\pi_2^{-1}(S)| + 2 g_D -2 d_2 g_C + 2d_2 -2 \text{ (by Hurwitz's formula)}\\
&\leq& |\pi_1^{-1}(S)| + 2 g_D -2 d_2 g_C + 2d_2 -2 \text{ (since $S$ is backward-complete)} \\
& \leq & |S| d_1 + 2 g_D -2 d_2 g_C + 2d_2 -2,
\end{eqnarray*}
from which the bound of the statement immediately follows.

Now if $S$ is a ramification-increasing backward-complete finite set,
\begin{eqnarray*} |S| d_2 &=& \sum_{z \in \pi_2^{-1}(S)} e_{2,z} \\
& \leq &  \sum_{z \in \pi_1^{-1}(S)} e_{2,z} \text{ (since $S$ is backward-complete)}\\
& \leq &  \sum_{z \in \pi_1^{-1}(S)} e_{1,z} \text{ (since $S$ is ramification-decreasing)}\\
&=& |S| d_1  \text{ (by formula~(\ref{sumS1}))}
\end{eqnarray*}
which under our assumption $d_2>d_1$ implies $S =\emptyset$.
\end{pf}
We record for later use a consequence of the proof:
\begin{scholium} \label{scholiumbalanced} Let $(D,\pi_1,\pi_2)$ be a self-correspondence with $d_1 \leq d_2$. If $S$ is a  ramification-increasing backward-complete finite set for $D$, then $S$ is complete and equiramified. 
Moreover, if a self-correspondence admits a complete equiramified non-empty finite set, then it is is balanced.
\end{scholium} 
 \begin{pf} The second chain of inequalities in the proof of the proposition is an equality by assumption $d_1 \leq d_2$, hence all intermediate inequalities must be equalities. The rest is clear.
 \end{pf}
 Applying the proposition (resp. the scholium) to ${}^t D$, we get a dual statement concerning forward-complete sets in the cases $d_1>d_2$ (resp. $d_1=d_2$) that we let the reader make explicit.
 Combining the proposition and its dual statement, we get:
\begin{cor} Let $(D,\pi_1,\pi_2)$ be a self-correspondence over $C$ of bi-degree $(d_1,d_2)$ with $d_1 \neq d_2$ (that is, $D$ is unbalanced). Let $d=max(d_1,d_2)$. Then any finite complete set $S$ is not equiramified, in particular is not étale, and satisfies
 $|S| \leq 2 \frac{g_D - d g_C +d -1}{|d_2-d_1|}$.
 \end{cor}

\begin{example} \label{examplegraphDf} Let $f \in k(X)$ be a rational function
of degree $d>1$, seen as a morphism $\P^1_k \rightarrow \P^1_k$. For concreteness and convenience we now explain and prove the classical results about the structure of the graph $\Gamma_{D_f}$ (those results are well-known to specialists and their proofs are probably all in the literature but  scattered in many different places). Note that $\Gamma_{D_f}$ is  a directed graph with at most $d$ entering edges and {\bf exactly one 
exiting edge} at each vertex.

\par \medskip
First introduce some more terminology concerning directed graphs.

A {\it rooted regular directed tree of valence $d+1$} is a directed graph with a specified vertex $r$ (the {\it root}), $d$ vertices $x_1,\dots,x_d$
with exactly one edge from $x_i$ to $r$, $d^2$ vertices $x_{i,j}$ ($i,j \in \{1,\dots,d\}$) with exactly one edge from $x_{i,j}$ to $x_i$,
$d^3$ vertices $x_{i,j,k}$ with exactly one edge from $x_{i,j,k}$ to $x_{i,j}$, and so on... 

An {\it unrooted regular directed tree of valence $d+1$} is
a directed graph such that the smallest backward-complete subset containing any vertex $r$ is a rooted directed tree of valence $d+1$ and root $r$. 

A {\it volcano} of type $(n,d)$ is a directed graph consisting of a cycle $C$ of length $n \geq 1$ and for every point $x \in C$, $d-1$ disjoint rooted directed trees of valence $d-1$ rooted at $d-1$ points $r_1,\dots,r_{d-1}$ which are exactly the backward-neighrbors of $x$ not in $C$.
 
 \par \medskip
Let $S$ be a complete irreducible set for $\Gamma_{D_f}$. We claim that
\begin{itemize}
\item[(i)] $S$ contains at most one cycle.
\item[(ii)] If $S$ is finite, then it is a cycle and every edge $z$ of $S$ is totally ramified (in particular, not étale). 
\item[(iii)] If $S$ is étale, then it is either an unrooted directed tree of valence $d+1$ or a volcano of type $(n,d)$ for some $n \geq 1$.
\end{itemize} 
To prove (i), if  there are two cycles $C$ and $C'$, let
$(x_0,\dots,x_n)$ be a shortest undirected path relying $C$ and $C'$. There must be an edge from $x_1$ to $x_0$ because there is already an exiting edge at $x_0$ in the cycle $C$. Similarly, there must be an edge from $x_2$ to $x_1$ because there is already an edge exiting at $x_1$ (and toward $x_0$), and similarly, for $(x_3,x_2), \dots, (x_n,x_{n-1})$. But there must be an edge starting from $x_n$ inside the cycle $C'$, a contradiction.

(ii) If $S$ is finite, since it is complete, it must contain a cycle $(x_0,\dots,x_n)$, and this cycle is unique by (i). If $S$ is not reduced to that cycle, there is a vertex $y_0 \in \Gamma_{D_f}$ not in the cycle with an edge from $y_0$ to some $x_i$ in the cycle, and $y_0 \in S$ by (backward-)completeness of $S$. There must be also a vertex $y_1 \in S$ with an edge from $y_1$ to $y_0$, and $y_1$ cannot be $y_0$ nor in the cycle because the unique edges starting at those points ends in the cycle. Then there must be a $y_2$ in $S$ with an edge to $y_1$, and $y_2$ can no be $y_1$, $y_0$ or in the cycle because the edges starting at those points ends in $y_1$ or in the cycle. We continue indefinitely in. the same way, to construct an infinite injective sequence of points $y_n$ in $S$, contradicting the finiteness of $S$.
Thus $S$ is a cycle. Since there are only one edge $z$ ending at each vertex of a cycle, this edge must be totally ramified: $e_{2,z}=d$.

For (iii), if $S$ is an étale complete set, and $x_0$ in $S$, then the $n$-backward-neighbors of $S$ are the $d^n$ distinct solutions of $f^{(n)}(x)=x_0$ (they are distinct, because the solutions are all in $S$, and for any $x \in S$, $(f^{(n)})'(x)=f'(x) f'(f(x)) \dots f'(x_0) \neq 0$ since $S$ is étale.). If an $n$-backward-neighbor $y_0$ of $x_0$ is the same as an $m$-backward neighbor of $z_0$ for some $n > m$, then $f^{(n-m)}(y_0)=y_0$ so $y_0$ belongs to a finite directed cycle of $S$, and so does $x_0=f^{(n)}(y_0)$.

Thus if $S$ has no cycle, the set of all $k$-backward-neighbors of any vertex $x_0 \in S$ is a rooted regular directed tree of valence $d+1$. It is then easy to see that $S$ itself is an unrooted directed tree of valence $d+1$.

Assume that $S$ has a finite directed cycle $C$, $x_0,\dots,x_n=x_0$. If $y_0 \in S$, then some $f^{(k)}(y_0) \in C$, otherwise $S$ would not be irreducible. Let $k \geq 0$ be the smallest such integer, and define a map $\pi: S \rightarrow C$ by $\pi(y_0)=f^{(k-1)}(y_0)$,
and another map $\pi':S-C \rightarrow S-C$, by $\pi'(y_0)=f^{(k-1)}(y_0)$ (the map is well-defined since one has $k \geq 1$ if $y_0 \not \in C$, and $f^{(k-1)}(y_0) \not \in C$ be minimality of $k$).
Then $\pi(y_0)$ is the forward-neighbor of $\pi'(y_0)$.  Let us determine the fibers of $\pi$, for instance $\pi^{-1}(x_0)$.
If $y_0 \not \in C$ ans $\pi(y_0)=x_0$, then $\pi'(y_0)$ is not in $C$, so it has to be one of the $d-1$ backward-neighbors of $x_0$ not in $C$, that we'll call $r_1,\dots,r_{d-1}$. Thus $\pi^{-1}(x_0) = \{x_0\} \cup \cup_{i=1}^{d-1} \pi'^{-1}(r_i)$.
The sub-graph $\pi'^{-1}(r_i)$ has no cycle, and is composed of $d^k$ $k$-backward-neighbors of $r_i$ for any $k \geq 0$. It is therefore a rooted directed tree of root $r_i$ and valence $d+1$. It follows that $S$ is a volcano of type $(n,d)$.
\par \medskip
Thus the complete irreducible sets are of 3 types: non-étales sets, étale volcanoes (of type $(n,d)$ for some $n$), and étale unrooted trees (of valence $(d+1)$).\begin{itemize} 
\item[(iv)] The number of complete irreducible sets that are non étale is finite, but non-zero. The number of  étale complete irreducible sets that are volcanoes is infinite countable; and the number of étale complete irreducible set that are unrooted directed trees is $0$ if and only if $k$ is not algebraic over
a finite field.
\end{itemize}
Indeed, since there are finitely many, but non-zero (since $d \geq 2$) non-étale points $z \in k$, the first assertion is clear. 

Note that a complete irreducible set contain a cycle of length $n$ if and only if it contains an $x \in \bar k$ such that $f^{(n)}(x)=x$.  Since $f^{(n)}$ has degree $d^n >1$, we may assume up to conjugating that $f^{(n)}$ does not fix $\infty$, hence is of the form $p(X)/q(X)$ where $p,q$ are polynomial with $\deg p < \deg q =d^n$. Hence $f^{(n)}(x) = x$ is equivalent to $p(x)-x q(x) = 0$, an equation of degree $d^n+1$ which always has a solution in $\bar k$. It follows $C(\bar k)$ always contains cycles of any length $n$, and finitely many of them for a given $n$, hence an infinite countable number of them.  By (i), there exists countably infinitely many irreducible complete sets containing a cycle, and since all of them but a finite number are étale, there are countably infinitely many étale volcanos.

For the assertion concerning unrooted trees, observe that an étale irreducible complete set is an étale unrooted tree if and only if it does not contain a directed cycle. If $k$ is algebraic over a finite field, and $x \in C(\bar k)$ then clearly $f$ and $x$ are defined over a finite field $k_0$, and since $f$ stabilizes the finite set $C(k_0)$, $x$ belongs to a directed cycle. It remains to show that when $k$ is not algebraic over a finite field, there are points that do not belong to a directed cycle. This will be done below in Prop.~\ref{proppolarized}.

\par \medskip

We can now describe the so-called {\it exceptional set} $E$:
\begin{itemize}
\item[(v)] The union $E$ of all complete finite sets has cardinality $0$, $1$ or $2$. If $E$ is a singleton, then we can make a change of variable moving $E$ to $\{\infty\}$, and then $f$ is polynomial. If $E$ is a pair, then we can make a change of variable moving $E$ to $\{0,\infty\}$, and then $f(z)=a z^{d}$ or $f(z)=az^{-d}$ for some $a \in k^\ast$. In the first case, $D_f$ has two finite irreducible complete sets,
$\{0\}$ and $\{\infty\}$, and in the second case one finite irreducible complete set, the pair $\{0,\infty\}$. 
\end{itemize}
Indeed the first assertion follows from Prop.~\ref{unbalancedprop} and the rest easily follows. (For another proof of (v) in the case $k=\C$ using complex analysis, see \cite[Theorem 3.6]{Milnor}.)
\end{example}

\subsection{Finitary self-correspondences}

\begin{theorem} \label{finitary}
For a self-correspondence $D$ on a curve $C$ over a field $k$, the following are equivalent: 
\begin{itemize} \item[(i)] There exists a non-constant $k$-morphism $h: C \rightarrow \P^1_k$ such that $h \circ \pi_1 = h \circ \pi_2$.
 \item[(ii)] There exists a non-constant $\bar k$-morphism $h: C_{\bar k} \rightarrow \P^1_{\bar k}$ such that $h \circ \pi_1 = h \circ \pi_2$
\item[(iii)] There is an integer $M$ such that every irreducible complete set has cardinality less or equal than $M$.
\item[(iv)] All irreducible complete sets of $D$ are finite.
\item[(v)] All irreducible complete sets of $D$ but possibly finitely many are finite.
\item[(vi)] There are infinitely many finite complete sets.
\end{itemize},
\end{theorem}
\begin{pf}
When $D$ is unbalanced, (i) and (ii) are false by additivity of the degree, and (iii), (iv), (v) and (vi) are false by Prop.~\ref{unbalancedprop}. The assertions (i) to (vi) are thus equivalent. We may therefore assume that $D$ is balanced of degree $d$.
\par \bigskip
The implication (i) $\Rightarrow$ (ii) is trivial. For (ii) implies (iii), we just note that the fibers $h^{-1}(t)$, $t \in \P^1(\bar k)$ are complete since $\pi_1^{-1}(h^{-1}(t))=(h \circ \pi_1)^{-1}(t)=(h \circ \pi_2)^{-1}(t)=\pi_2^{-1}(h^{-1}(t))$, and they all have size $\leq \deg h$. Therefore every irreducible complete set has  cardinality less or equal than $\deg h$, which gives (iii).

That (iii) implies (iv) and (iv) implies (v) is trivial; and (v) implies (vi) by Prop.~\ref{infinitecomplete}.

It is also not hard to prove that (ii) implies (i), as in \cite[Prop. 3.8]{raju}. In fact, assuming (ii), we know that there exists a finite extension $k'$ of $k$ and a map $h'$ from $C_{k'}$ to $\P^1_{k'}$ defined over $k'$ such that $h' \circ \pi_1 = h' \circ \pi_2$. Let $V$ be the Weil's restriction of scalars of $\P^1_{k'}$ to $k$, and $\tilde h: C \rightarrow V$ the $k$-morphism corresponding to $h$ according to the universal property of Weil's restriction. The image $C'$ of $C$ by $\tilde h \circ \pi^1=\tilde h \circ \pi_2$ in $V$, with its reduced scheme structure, is a closed subscheme of $V$ defined over $k$ and  which has positive Krull dimension. There is therefore a rational function $u \in k(V)$ that induces a non-constant map on $C'$. Thus, if $h := u \circ \tilde h$, $h$ is a $k$-morphism $C \rightarrow \P^1_k$ such that $h \circ \pi_1 = h \circ \pi_2$ and (i) holds.

\par \bigskip
It only remains to prove that (vi) implies (ii).

 To prove this, we first reduce to the case where $k$ is of finite type over its prime field. There is a subfield $k_0$ of $k$, of finite type over the prime subfield of $k$ such that $C$, $D$, $\pi_1$ and $\pi_2$ are defined over $k_0$. Assuming (vi), that is that $D(\bar k)$ contains infinitely many finite irreducible  complete sets, we must show that (vi) holds for $k_0$, that is $D(\bar k_0)$ also contains infinitely many finite irreducible  complete sets. If all the finite complete sets in $D(\bar k)$ are in $D(\bar k_0)$, we are done. Otherwise, there is one finite complete set $S$ in $D(\bar k)$ which is not in $D(\bar k_0)$. Let $k'$ be a subfield of $\bar k$ containing $k_0$ and of finite type over $k_0$ such that all points of $S$ are defined over $k'$. By assumption, $k'$ is transcendental over $k_0$. Let $V=\spec A$ be an affine algebraic variety over $k_0$ whose field of fraction is $k'$ and such that $S$ spreads out to $V$, that is that all points of $S$ are defined over $A$. Then $V$  has positive Krull's dimension, so it has infinitely many $\bar k_0$-points. Any $\bar k_0$-point $t$ in $V$ gives by specialization of $S$ a finite complete set $S_t$ in $D(\bar k_0)$, and the $S_t$ are all distinct. It follows that $D(\bar k_0)$ contains infinitely many finite complete sets, and (vi) holds for $k_0$. Replacing $k$ by $k_0$, we may henceforth assume that $k$ is of finite type over  its prime field.

We now complete the proof of (vi) $\Rightarrow$ (ii) assuming $k$ is finitely generated over the prime field. Let $k^\sep$ and $k^\perf$ be the separable and perfect closures of $k$ in $\bar k$. Then $\bar k = k^\sep \otimes_k k^\perf$ and $\Gal(k^\sep/k)=\Gal(\bar k/k^\perf)$. We denote this Galois group by $G$.

Let $J$ be the Jacobian of $C$ over $k$. Since $k$ is of finite type over its prime field, $J(k)$ is a finitely generated abelian group
by the theorem of Néron (\cite{Neron}, see also \cite{LN}). One has $J(k^\perf) \otimes_\Z \Q = J(k) \otimes_\Z \Q$. (Indeed, there is nothing to prove in characteristic $0$, and in characteristic $p>0$, $J(k^\perf) = \cup_n J(k^{1/p^n})$, so by induction it suffices to prove that $pJ(k^{1/p}) \subset J(k)$. But $$p J(k^{1/p}) = V F  J(k^{1/p})  \subset V J^{(p)}(k) \subset J(k)$$ where $F: J \rightarrow J^{(p)}$ and $V: J^{(p)} \rightarrow K$ are the relative Frobenius and Verschiebung maps.) Let $r$ be the finite dimension of $J(k^\perf) \otimes_\Z \Q$.

By assumption, there are infinitely many finite complete sets in $C(\bar k)$, and therefore infinitely many étale finite  complete sets. Each of them is a finite subset of $C(\bar k)$, hence has a finite $G$-orbit. By grouping the irreducible complete sets by $G$-orbits, we see that there are still infinitely many disjoint étale finite complete sets invariant by $G$. Let us chose $r+2$ of them, say $S_0,\dots,S_{r+1}$.
To every $S_i$ we attach the Weil divisor $$\Delta_i = \sum_{x \in S_i} [x]\in \Div \, C_{\bar k}$$ and let $\delta_i = |S_i|$ be its degree. The $r+1$ divisors $\delta_0 \Delta_i - \delta_i \Delta_0$ have degree zero, hence they define points in $\Pic^0(C_{\bar k})=J(\bar k)$, where $J$ is the Jacobian of $C$, and those points are $G$-invariant, hence in $J(k^\perf)$. Therefore those points are $\Q$-linearly dependent, hence $\Z$-linearly dependent, which means that there are non-all-zero integers $n_i$, $i=0,\dots,r+1$, and a non-constant  $k^\perf$-map $h: C_{k^\perf} \rightarrow \P^1_{k^\perf}$ such that $$\sum_{i=0}^{r+1} n_i \Delta_i  = \div\, h.$$ 
Since the $S_i$ are étale,  for $j=1,2$ one has $$\pi_j^\ast \Delta_i = \sum_{z \in \pi_j^{-1}(S_s)} [z],$$ and since the $S_i$ are complete,
it follows that $$\pi_1^\ast \Delta_i = \pi_2^\ast \Delta_i.$$ It follows that $$\div (h \circ \pi_1)=\pi_1^\ast \div\, h = \pi_2^\ast \div\, h  = \div (h \circ \pi_2),$$ hence that there exists $\lambda \in (k^\perf)^\ast$ such that $$\lambda h \circ \pi_1 = h \circ \pi_2.$$ This implies that (reasoning as in the proof of (ii) implies (iii)) if $S$ is any complete set, $h(S)$ is stable by multiplication by $\lambda$ and $\lambda^{-1}$. By assumption, there exists a finite complete $S$ such that $h(S)$ is not contained in $\{0, \infty\}$. This implies that $\lambda$ is a root of unity, and if $\lambda^n=1$, replacing $h$ by $h^n$ gives (ii). This completes the proof of the implication (vi) implies (ii), hence of the theorem.
\end{pf}

Remember from the introduction that when the equivalent assertions of the preceding theorem are satisfied, we say that $D$ is 
{\it finitary}.

\begin{remark} A self-correspondence $D$ is finitary if and only if its transpose ${}^t D$ is finitary. This is seen trivially on any of the assertion (i) to (vi).

Also, if $k'$ is any field containing $k$, a self-correspondence $D$ is finitary if and only if its base change $D_{k'}$ is finitary. The implication $D$ finitary $\Rightarrow D_{k'}$ finitary is clear on (i), and its converse is clear on (iii) or on (iv). 
\end{remark}
\begin{remark} A correspondence of morphism type, say $D_f$, is finitary if and only if $f$ is an automorphism of finite order. Indeed, if $D_f$ is finitary, then $\deg f=1$ by (i) and $f$ is an automorphism, whose action on the generic fiber of $h$ is a bijection of a finite set, so some power of $f$ acts trivially on the generic fiber of $h$, hence on $C$, and $f$ has finite order. The converse is trivial.
\end{remark}

\begin{remark} \label{remarkraju} The method of using Jacobians in the proof of the Theorem is inspired by \cite[chapter 9]{raju}, and our condition (a) is inspired by the condition he calls ``having a core"\footnote{Raju himself is inspired by Mochizuki, which introduces the same notion in the hyperbolic case, under the name ``having an hyperbolic core", see \cite{mochizuki}. In the case $k=\C$, a closely related notion has also been considered by Bullett, Penrose, and their coauthors (see e.g. \cite[\S2.5]{BP}), under the name of {\it separable} (self-)correspondence. This notion is equivalent to ``having a core" in the minimal case. Apparently, the two sets of authors (Bullett et al., Mochizuki and Raju) were unaware of each other's works.}
 which gives the title to his article. For Raju, a (general, not self-) correspondence $(D,\pi_1:D\rightarrow C,\pi_2:D \rightarrow C')$ {\it has a core} if there exist non-constant maps $f:C \rightarrow \P^1$, $g: C' \rightarrow \P^1$ such that $f \circ \pi_1=g \circ \pi_2$. For a self-correspondence, his notion is much weaker than our notion of being finitary, where we require $f=g$. Indeed, finitary implies balanced, while there are plenty of unbalanced self-correspondences that have a core, in particular all those of morphism type. Even for balanced self-correspondences, having a core does not imply being finitary, as the example of $D_f$ when $f$ is an automorphism of infinite order shows (or for less trivial examples, a suitable étale self-correspondence on a curve $C$ of genus $1$, for such a correspondence always have a core in the sense of Raju, but in general is not finitary). In the case of a symmetric self-correspondence $D$, however, one can show that $D$ is finitary if and only if it has a core in the sense of Raju: see Lemma~\ref{lemmacorefinitary}. 
\end{remark}

\begin{remark} \label{remarkDKW} One finds in the literature yet another property akin to ``having a core" or ``being finitary", namely what 
Dinh, Kaufmann et Wu calls ``weakly modular" in \cite{DKW}. They work in the case $k=\C$ and they say that a balanced self-correspondence $D$ over a curve $C$ is {\it weakly modular} if there are two probability measures $m_1$ and $m_2$ over $C(\C)$ such that $\pi_1^\ast m_1 = \pi_2^\ast m_2$. To ``complete
the square", let us say that $D$ is {\it weakly finitary} if there is one probability measure $m$ on $C(\C)$ such that $\pi_1^\ast m = \pi_2^\ast m$. Thus one has a following square of implications for $D$ a balanced self-correspondence over $\C$, none of which being an equivalence.
$$\xymatrix{ D \text{ is finitary} \ar@{=>}[r] \ar@{=>}[d] &  D \text{ has a core} \ar@{=>}[d] \\ D \text{ is weakly finitary} \ar@{=>}[r] & D \text{ is weakly modular} } $$
Note that any self-correspondence which has a non-empty finite complete set  $S$ is weakly finitary, hence weakly modular, as taking $m$ the normalized counting measure on $S$ shows. Thus to be weakly modular is a very weak condition.

Theorem 1.1 of \cite{DKW} states that for any balanced non-weakly-modular self-correspondence $D$ on a  curve $C$ over $\C$, there exists a measure $\mu_D$ on $C(\C)$ which does not charge polar sets (in particular finite sets), and such that for every smooth measure\footnote{We recall the standard measure-theoretic notation used here: let $\pi: D(\C) \rightarrow C(\C)$ be an holomorphic map, non-constant on every component of $D$; if $\mu$ is a Borel measure on $D(\C)$, then $\pi_\ast \mu$ is the measure on $C(\C)$ defined by $\pi_\ast(\mu)(B)= \mu(\pi^{-1}(B))$; if $\mu$ is a Borel measure on $C(\C)$, then $\pi^\ast \mu$ is the Borel measure on $D(\C)$ defined by $\int f d\pi^\ast \mu = \int \pi_\ast f d \mu$ for every continuous function $f$ on $D(\C)$, where $\pi_\ast f(x) = \sum_{z \in \pi^{-1}(x)} f(x)$.}
$\mu$ on $C(\C)$, $(\frac{1}{d} (\pi_1)_\ast (\pi_2)^\ast)^n \mu \rightarrow \mu_D$ as $n \rightarrow \infty$. In the case of a non-balanced self-correspondence with $d_1<d_2$, the theorem was known earlier (\cite{BS}). One may ask whether this theorem  holds more generally for any $D$ that is not finitary.  We conjecture the answer to be
 yes. (See also Remark~\ref{remarkDKW2}.)
\end{remark}

\subsection{The number of irreducible complete finite sets, I} \label{finiteetaleI}

Given a  balanced non-finitary self-correspondence $D$ over $C$ of bi-degree $(d,d)$, can we give a bound to the number of irreducible complete finite sets in terms of the genera $g_C$ and $g_D$ of the curve involved and the degree $d$?  The proof of Theorem~\ref{finitary} does not provide such a bound, even involving not only $g_C$, $g_D$, and $g$, but the Mordell-Weil rank  $r$ of the Jacobian variety of $C$ over $k_0$, because the step where we group together finite complete sets to get $G$-invariant
complete sets is not effective.

This question can be separated in two, one concerning the number of {\it étale}  irreducible complete finite sets and one about the non-étale ones.

Concerning the number of non-étale irreducible complete finite sets, it is clearly bounded by the number of ramification points of $\pi_1$ and $\pi_2$, which in turns is bounded, using Hurwitz's formula, by $2g_D-2-d(2g_C-2)$. 

The question  of bounding the number of étale irreducible complete finite sets is more subtle. In this section, we give several such bounds
in particular but important situations (namely the case where $k$ is a finite field, or when $C=\P^1$, or when the correspondence $D$ is symmetric) using effective variants of the proof of Theorem~\ref{finitary}.
In the next section (\S\ref{finiteetaleII}), using completely different methods we give a general result in characteristic zero (and also under a weaker form in characteristic $p$), namely that a non-finitary correspondence has at most $2$ étale (or even equiramified) finite irreducible sets.

\begin{prop} \label{propetalefinite} Assume that $k$ is algebraic over a finite field.
If a correspondence $D$ on a curve $C$ over $k$ is not finitary, then it has at most one non-empty finite equiramified complete set, and in particular at most one étale non-empty complete set.\end{prop}
Note that the proposition implies that in case there is one non-empty equiramified complete set, it is automatically irreducible. Of course, $D$ can very well have no non-empty finite complete set, as in the case of $D_f$, where $f$ is a the translation by a non-torsion element on an elliptic curve.

\begin{remark} The hypothesis made on $k$ is necessary. Indeed, assume that $k$ is not algebraic over a finite field. 
Then $k$ has an element $t \neq 0$ which is of infinite multiplicative order (take an element which is transcendental over $\F_p$ if $k$ has characteristic $p$ and $t=2$ if $k$ has characteristic 0). Consider the map $f(x)=tx$ from $\P^1$ to $\P^1$, and the correspondence $D_f$ it defines. This correspondence is of bi-degree $(1,1)$, and has two complete finite sets, obviously étale: $\{0\}$ and $\{\infty\}$. Yet $D_f$ is not finitary, because $f$ is not an automorphism of finite order.
\end{remark}

\begin{pf}
Suppose that there are two distinct finite equiramified non-empty complete sets $S$ and $S'$. If $S \subset S'$, replace $S'$ by $S'-S$. If $S \not \subset S'$, replace $S$ by $S-(S'\cap S)$. This way we may assume that $S$ and $S'$ are not only distinct, but disjoint.

In view of our hypothesis on $k$, we may assume that the self-correspondence is defined over a finite field $k_0$,  that $S$ and $S'$ are subsets of $C(k_0)$, and that $\pi_1^{-1}(S)$ and $\pi_2^{-1}(S')$ are subsets of $D(k_0)$.

Let  $\Delta_S = \sum_{s \in S} [s]$ and $\Delta_S' = \sum_{s \in S' }[s]$ be the effective Weil's divisors attached to $S$ and $S'$. The divisor
$|S'| \Delta_S - |S| \Delta_{S'}$ has degree zero, hence is torsion in $\Pic_{k_0}(C)$ (since it is an element of the finite group $\Pic^0_{k_0}(C)$, the group of $k_0$-rational points of the Jacobian of $C$). Therefore, there exist $n$ and $m$ such that $n \Delta_S - m \Delta_{S'}$ is a principal divisor; in other words, there exists a rational function $h$ in $k(C)$ such that $\div\,  h= n\Delta_S - m \Delta_{S'}$.

We claim that there exists $\lambda \in k^\ast$ such that $h \circ \pi_1 = \lambda h\circ \pi_2$ as functions in $k(D)$, up to multiplication by a non-zero scalar. Indeed, for $i=1,2$, $\div\, h \circ \pi_i = \pi_i^\ast (n \Delta_S - m \Delta_{S'}) = \sum_{t \in \pi_i^{-1}(S)} n e_{i,t }[t] - \sum_{t \in \pi_i^{-1}(S')} m e_{i,t'} [t']$ where $e_{i,t}$ is the ramification index of $\pi_i$ at $t$. Using that $S$ and $S'$ are equiramified complete sets, we see that $\div\, h \circ \pi_1 = \div\, h \circ \pi_2$, hence the claim.

Now since $\lambda$ belongs to a finite field, there is an $n$ such that $\lambda^n=1$. Replacing $h$ by $h^n$, we see that $D$ satisfies assertion (i) of Theorem~\ref{finitary}.
\end{pf} 
 
\begin{prop} \label{caseP1} Let $k$ be any field. If $D$ a self-correspondence over $\P^1_k$ which is not finitary, there are at most two irreducible finite equiramified complete sets, and when there are two of such, there is no other irreducible finite complete set at all.
\end{prop}
\begin{pf}
Arguing as in the case of a finite base field $k$,  but using the fact that $\Pic^0 C(k)$ is trivial, we see that if $S$ and $S'$ are two disjoint irreducible equiramified sets,  there is a function $f : C \rightarrow \P^1$ with divisor $|S'| \Delta_S-|S|\Delta_{S'}$. Such a function satisfies $f \circ \pi_1 = \lambda f \circ \pi_2$ for some $\lambda \in k^\ast$, and also by definition $f^{-1}(\infty)=S$, $f^{-1}(0)=S'$. Since $D$ is not finitary, $\lambda$ is not a root of unity, and as in the proof of Theorem~\ref{finitary},
we see that for $T$ finite complete, $f(T) \subset \{0,1\}$. Since  $f^{-1}(0)=S$ and $f^{-1}(\infty)=S'$ are irreducible, they are the two only irreducible complete sets.
\end{pf}

Finally, we give two more results in the same vein but for symmetric self-correspondence. They are based on earlier results in the literature and the following lemma.

\begin{lemma} \label{lemmacorefinitary} If a self-correspondence has a core in the sense of Raju (see Remark~\ref{remarkraju}, or \cite[Definition 3.5]{raju}) then ${}^t D D$ is finitary. If moreover $D$ is symmetric then it is itself finitary. \end{lemma}
\begin{pf} If $(D,\pi_1,\pi_2)$ has a core, that is if there exists $f,g: C\rightarrow \P^1$ such that $f \circ \pi_1 = g \circ \pi_2$, then the forward map of $D$ sends fibers of $f$ to fibers of $g$, and the backward map of $D$, i.e. the forward map of ${}^t D$, sends fibers of $g$ to fibers of $f$. Thus ${}^t D D$ preserves the fibers of $f$, and is therefore finitary, which proves the first assertion.

If moreover $D$ is symmetric, then $D^2$ is finitary and has infinitely many irreducible complete finite sets. But every irreducible complete set $S$ of $D$ breaks down in at most two  irreducible complete sets of $D^2$,  namely the set of points of $S$ which are connected to a given point $x_0$ of $S$ by a path of even length, and its complement if non-empty.
Therefore $D$ must have infinitely many finite complete sets, and is therefore finitary.
 \end{pf}

\begin{prop}[Raju] \label{propraju} If $D$ is a non-finitary symmetric self-correspondence on a curve $C$ over a field $k$, then $D$ has at most one irreducible finite equiramified complete set. If $k$ has characteristic zero and $\pi_1$, $\pi_2$ are étale, then $D$ has no  irreducible finite complete set.
\end{prop}
\begin{pf} Raju proves that if $D$ is a self-correspondence without a core, it has at most one finite irreducible finite étale complete set (see \cite[Theorem 9.6]{raju}) and none in characteristic $0$ when $\pi_1,\ \pi_2$ are étale (see \cite[Cor 9.2]{raju}). In fact his proofs work with ``étale complete'' replaced with ``equiramified complete'', as in the proofs of Prop.~\ref{propetalefinite} and Prop.~\ref{caseP1}.  By the lemma above, this implies the proposition.
\end{pf}
\begin{cor} Let $D$ be a non-finitary symmetric self-correspondence on a curve $C$ over a field $k$. Let $S$ be an irreducible finite complete étale set for $D$. Then the undirected graph of $S$ (obtained by forgetting the orientation of the edges) is not bipartite.
\end{cor}
\begin{pf} If the undirected graph $\Gamma_S$ attached to $S$ was bipartite, then the graph $\Gamma_{S,2}$ with the same set of vertices $S$ but whose edges are paths of degree $2$ in $S$ would be disconnected. But $D^2$ is also not finitary (by Lemma~\ref{lemmacorefinitary}), and its graph on the set of vertices $S$ is $\Gamma_{S,2}$ (by Lemma~\ref{lemmaSDD'}, using that $S$ is étale), in contradiction with Prop.~\ref{propraju}.
\end{pf}

\begin{remark} Let $S \in C(\bar k)$ be a finite complete set for a self-correspondence $D$ on $C$. 
We shall say that the set $S$ is {\it consistently ramified} if for every undirected cycle with edges $z_1,\dots,z_n$, the rational number $\prod_{i=1}^n \left( \frac{e_{2,z_i}}{e_{1,z_i}}\right)^{\epsilon_i}$ is $1$, where the signs $\epsilon_i$ (for $i=1,\dots,n$) are defined to be $+1$ is the edge $z_i$ has a compatible orientation with $z_{i+1}$ (i.e. the target of $z_i$ is the source of $z_{i+1}$ or the source of $z_i$ is the target of $z_{i+1}$)  or $-1$ otherwise (we use the convention that $z_{n+1}=z_1$).
If $S$ is equiramified, it is consistently ramified, since all factors in the above product are 1; but clearly the converse is false.

We claim that Propositions~\ref{propetalefinite},~\ref{caseP1}, and~\ref{propraju} are still true with the phrase ``finite equiramified complete sets" replaced with ``finite consistently ramified complete sets". Indeed, if $S$ is consistently ramified, it is easy to see that one can attach to every $s \in S$ a positive integer $n_s$ such that for every $z \in \pi_1^{-1}(S)=\pi_2^{-1}(S)$, one has $n_{\pi_1(z)} e_{1,z} = n_{\pi_2(z)} e_{2,z}$. In the proof of prop.~\ref{propetalefinite} (for example), it suffices to change the definition of $\Delta_{S}$ to be the divisor $\sum_{s \in S} n_s [s]$, to obtain a divisor with support $S$ and such that $\pi_1^\ast \Delta_S=\pi_2^\ast \Delta_S$, and the rest of the proof may remain unchanged.
\end{remark}

A self-correspondence is {\it critically finite} (see \cite{B1}) if every ramification point in $C(\bar k)$ of $\pi_1$ or $\pi_2$ 
belongs to a finite complete set. In other words, the union $E_\crit$ of the irreducible non-étale complete sets is finite.

\begin{prop} \label{mochi} Let $k$ be a field of characteristic zero. Let  $D$ be a symmetric  non-finitary critically finite self-correspondence on $C$ over $k$. Assume that the curve $D$ is irreducible. Then $D$ has no non-empty étale finite complete set.
\end{prop} 
\begin{pf} (Inspired by section 3 of \cite{mochizuki}.)

We may assume $k=\bar k=\C$. If $E_\crit$ is empty, then $\pi_1,\pi_2$ are étale, and the result follows from the Prop.~\ref{propraju}.
Let $C_0 := C-E_\crit$, and $D_0 = D - \pi_1^{-1}(E_\crit)=D-\pi_2^{-1}(E_\crit)$, and we still denote by $\pi_1$ and $\pi_2$ the restriction of $\pi_1$ and $\pi_2$ to $D_0$. They are étale maps, and $(D_0,\pi_1,\pi_2)$ is, in an obvious sense, a self-correspondence over $C_0$ in the category of open curves.

If $C=\P^1$, and $E_\crit$ has $1$ or $2$ elements, then $C_0$ is the affine line or punctured affine line, and it has at most one étale finite cover of any degree $d$. Thus $\pi_1=\pi_2$ contradicting the assumption that $D$ is not finitary.

In the remaining cases, the open curve $C_0 := C-E_\crit$ is hyperbolic. Let $D_0 = D - \pi_1^{-1}(E_\crit)=D-\pi_2^{-1}(E_\crit)$. Let us identify the universal cover of $D_0$ (which is also a universal cover of $C_0$) with the upper half-plane  $\poincare$.

 Fix some $x \in C_0$, $z_1$ in $\pi_1^{-1}(x) \in D_0$  and $h_1 \in \poincare$ a point that maps to $z_1$ in $D_0$. The fundamental groups $\pi_1(C_0,x)$ and $\pi_1(D_0,z_1)$ are canonically identified, after those choices, with discrete subgroups of $\PSL_2(\R)$ that we shall denote respectively by $\Gamma_C$ and $\Gamma_D$; we have $\Gamma_D \subset \Gamma_C$, the inclusion being of finite index. Choose also a $z_2 \in \pi_2^{-1}(z) \in D_0$, $h_2 \in \poincare$ that maps to $z_2$, and a $g \in \PSL_2(\R)$ such that $g z_1=z_2$. Thus $\pi_1(D_0,z_2)$ is canonically identified with $g \Gamma_D g^{-1}$, which also is a subgroup of finite index of $\Gamma_C$, and we have a commutative diagram
 $$\xymatrix{ &  \poincare/\Gamma_D \ar[d] \ar[rrr]^{\sim}  & & & D_0 \ar[ddl]^{\pi_1} \ar[ddr]_{\pi_2} & \\
                 & \poincare/(\Gamma_C \cap  g^{-1} \Gamma_C g) \ar^{s_1}[dl] \ar_{s_2}[rd]& & & & \\
                \poincare/\Gamma_C \ar^{\sim}@/_1pc/[rrr] & & \poincare/\Gamma_C   \ar^{\sim}@/_1pc/[rrr] & C_0  & & C_0}$$
         where the unnamed horizontal maps (curved or straight)  are the identifications fixed above, and the diagonal maps $s_1$ and $s_2$ are given by respectively the inclusion and the conjugation by $g$ of  $\Gamma_C \cap  g^{-1} \Gamma_C g$ into $\Gamma_C$,
                the vertical map being given by the inclusion $\Gamma_D \subset \Gamma_C \cap g^{-1} \Gamma_C g$.
                This inclusion shows that $\Gamma_C \cap g^{-1} \Gamma_C g $ has finite index in $\Gamma_C$ and thus that $g$ belongs to the commensurator of $\Gamma_C$ in $\PSL_2(\R)$.

Let $\Gamma$ be the closure of the subgroup of $\PSL_2(\R)$ generated by $\Gamma_C$ and $g^{-1} \Gamma_C g$. We claim that $\Gamma$ has infinite index in $\Gamma_C$. Otherwise, $\Gamma$ would also be a lattice, and we would have two finite étale maps (surjective, of analytic stacks of dimension 1)  $ \poincare / \Gamma_C \rightarrow \poincare/\Gamma$ given by the inclusion and the conjugation by $g^{-1}$ of $\Gamma_C$ into $\Gamma$. We could then give an algebraic structure on $\poincare/\Gamma$, making it an algebraic stacks, and choose a non-constant map of $\poincare/\Gamma$ to $\P^1$, which composed with the two finite étale maps  $ \poincare / \Gamma_C \rightarrow \poincare/\Gamma$ gives two morphisms of Riemann surfaces $f,g : \poincare/\Gamma_C = C_0 \rightarrow \P^1$ such that $f \circ \pi_1 = g \circ \pi_2$. 
Extending $f,g$ to the complete smooth curves $C$, we thus see that  the symmetric correspondence $(D,\pi_1,\pi_2)$ has a core $(f,g)$, hence is finitary by Lemma~\ref{lemmacorefinitary}, contradicting our hypothesis.

We just need to show that the self-correspondence of Riemann surfaces $(D_0,\pi_0,\pi_1)$ has no finite complete sets,
and by the commutativity of the diagram above it suffices clearly to show hat the self-correspondence $\poincare/(\Gamma_C \cap g \Gamma_C g^{-1}, s_1,s_2)$ on $\poincare/\Gamma_C$ has no finite complete set. If $S$ was such a complete set, its preimage $\tilde S$ in $\poincare$ would be invariant by $\Gamma_D$ (obviously) and by $g \Gamma_D g^{-1}$, hence by $\Gamma$. The set $\tilde S$ would thus be an infinite union of $\Gamma_C$-orbits, contradicting the finiteness of $S$.
\end{pf}
 
\subsection{Complements and questions}

\subsubsection{On bounding the size of a finite complete set}	

Given a balanced non-finitary self-correspondence $D$ over $C$, the size of any finite complete set is bounded by the size of the union of all
finite complete sets, the finite {\it exceptional set} $E$. 
The question of bounding the size of $E$ seems much more difficult than bounding the number of irreducible complete finite sets (i.e. the number of components of $E$). The method of proof of Prop.~\ref{unbalancedprop} breaks down in the balanced case, and we do not know any general such bound.

We content ourselves by giving one result valable on $\P^1$ and in characteristic zero. It is a rephrasing of a  result due to Pakovitch.

\begin{prop} Let $k$ be a field of characteristic zero, $(D,\pi_1,\pi_2)$ a self-correspondence over $\P^1_k$ of bi-degree $(d,d)$. Let $g_D$ be the genus of $D$. We assume that 
\begin{itemize} \item[(i)] The singleton $\{\infty\}$ is a complete equiramified set.
\item[(ii)] There exits a $\lambda \in k^\ast$ such that for every $z \in \pi_1^{-1}(\infty)=\pi_2^{-1}(\infty)$, $$\ord_z (\pi_1 - \lambda \pi_2) > \ord_z \pi_1 = \ord_z \pi_2.$$
\item[(iii)] $D$ is not finitary.
\end{itemize}
Then every complete finite set has cardinality $ \leq 3 + (2g_D-1)/d$, with equality possible only in the case $g_D=0$, $d=1$.
\end{prop}
\begin{pf}
This is essentially Théorème 1 of \cite{pakovitch}.
More precisely, to prove our proposition we reduce to the case $k = \C$. Assuming that the proposition is false, there is a finite complete set in $\P^1_\C$ 
of cardinality $> 3+(2g_D-1)/d$, hence, removing $\infty$, there is a finite complete set $K \subset \C$ of cardinality $> 2+(2g_D-1)/d$. Conditions (i) and (ii) allow us to apply Theorem 1 of Pakovitch which tells us that there is a rotation $\sigma$ of the plane $\C$ such that $\sigma(K)=K$ and $\pi_1 =\sigma \circ \pi_2$. Since $K$ is finite, $\sigma^{(\# K) !}$ is the identity of $K$, hence of the real affine closure of $K$,
and since $\# K \geq 2$,  the rotation $\sigma^{(\# K) !}$ must fix a real line in $\C$, hence is the identity of $\C$. Thus $h(z)=z^{(\# K) !}$ satisfies $h \circ \pi_1 = h \circ \pi_2$, in contradiction with (iii). 
\end{pf}
\begin{remark} The condition (ii) is quite restrictive in practice. Assuming (i), it is clear that  for every $z \in \pi_1^{-1}(\infty)=\pi_2^{-1}(\infty)$, there exits a $\lambda \in k^\ast$ such that $\ord_z (\pi_1 - \lambda \pi_2) > \ord_z \pi_1 = \ord_z \pi_2,$ but the existence of such a $\lambda$ {\it independent of $z \in \pi_1^{-1}(\infty)$} is problematic -- except of course when $\pi_1^{-1}(\infty)$ is a singleton.

Moreover, the theorem is false without condition (ii), as the following example given by Pakovitch shows: $D=\P^1$, $\pi_1(z) = \frac{z^2-z-1}{z^2+z+1}$, $\pi_2(z)=-\frac{z^2+3z+1}{z^2+z+1}$. Then $\pi_1^{-1}(\{\infty\})=\pi_2^{-1}(\{\infty\})=\{j,j^2\}$ and $\{\infty\}$ is  complete étale so (i) is satisfied, and one can check that (iii) is also satisfied.
But $\{-1,1\}$ is also a finite complete set (since $\pi_1^{-1}(\{-1,1\}) = \{-1,0,\infty\} = \pi_2^{-1}(\{-1,1\})$). This provides an example where $C=D=\P^1$, and $|E|\geq 3$.
\end{remark}

%

\subsubsection{Backward exceptional kernels}

If $\Gamma=(V,E,s,t)$ is a directed graph, we define the {\it backward exceptional kernel} $\Kb$  as the union in $V$ of all finite backward-complete sets.
A symmetric definition could of course be given for {\it forward exceptional sets} and we will let the interested reader reformulate the results below in this case.

If $D$ is a self-correspondence over $C$ of bi-degree $(d_1,d_2)$, its {\it backward exceptional kernel}  is the one of its associated directed graph $\Gamma_D$. A natural question for a self-correspondence is then: when is $\Kb$  finite? We cannot offer a complete answer to this question. Here is the little the author knows on $\Kb$.

When $d_1<d_2$, Prop.~\ref{unbalancedprop} shows that $\Kb$ is finite and gives a bound to its size.

What about $\Kb$ when $d_1 \geq d_2$? 

Consider first the case $d_2=1$, that is of a transpose of self-correspondence of morphism type :  $D \simeq {}^t D_f$.
If $d_1=1$, then $f$ is an automorphism of $C$, of infinite order, and $\Kb$ is finite. If $d_2>1$ however, then $\Kb$ is always countable infinite. Indeed, in this case we are in the situation of Example~\ref{examplegraphDf}, and by assertions (iii) and (iv) of that example, $\Kb$ contains the union $V$ of all volcanoes of $f$, which is infinite countable and is contained in the union of $V$ with finite union of the non-étale irreducible complete sets (which is countable). 

Finally, in the case where $d_1 \geq d_2 > 1$, $\Kb$ may be finite:
for instance, consider the case $D=C=\P^1$, $\pi_1(x)=x^3$, $\pi_2(x)=x^2$ of bi-degree $(3,2)$, where it is easy to see that $\Kb=\{0,\infty\}$, or for a balanced case, any symmetric non-finitary correspondence (e.g. an Hecke correspondence $D_l$ on the Igusa curve), where by symmetry, $\Kb=E$ which is finite. One can also construct trivial examples with $\Kb$ infinite, for example the sum ${}^t D_f + {}^t D_f$, where $f$ is an endomorphism of degree $2$, which has bi-degree $(2,4)$ and, like ${}^t D_f$, an infinite countable $\Kb$.
However, this self-correspondence is not minimal. This suggests:
\begin{question} Does there exist a minimal non-finitary self-correspondence of bi-degree $(d_1,d_2)$ with $ d_1 \geq d_2 > 1$ and $\Kb$ infinite?
\end{question}
 To analyze this question, note that for $S$ a complete irreducible subset of $C(\bar k)$, $S \cap \Kb=\Kb(S)$, and $\Kb = \coprod_S \Kb(S)$ when $S$ run among irreducible complete subsets. If $S$ is equiramified, then $\Kb(S)$ is a union of equiramified finite backward-complete subsets of $S$, and those subsets are complete by Scholium.~\ref{scholiumbalanced}; thus, if $\Kb(S)$ is not empty it is finite and equal to $S$. This shows that the answer to the question is yes when $\Kb$ is equiramified, and because the number of non-equiramified irreducible complete sets $S$ is finite, in general the question reduces to ``is $\Kb(S)$ finite when $S$ is a non-equiramified complete set?"

\subsubsection{Backward exceptional and forward exceptional sets}

If $D$ is a self-correspondence over $C$ of bi-degree $(d_1,d_2)$, we define the {\it backward exceptional set} $\Eb$ as the  smallest forward-complete
set containing $\Kb$. Be careful that the {\it backward exceptional set} is forward-complete by definition, but not in general backward-complete.

If $D$ is not finitary, $\Eb$ is always ``small'': it contains the finite exceptional set $E$  and is contained in the union of $E$ and finitely many non-equiramified irreducible complete sets. In particular, $\Eb$ is  at most countable, and its complement contains an infinite union of irreducible complete sets. This follows easily from our study of $\Kb$.

If $\Kb$ is étale, then it is complete and  $\Eb=\Kb$ is finite. However, $\Eb$ may be infinite in general. An example due to Dinh and Favre showing this is given in \cite[Exemples 3.11]{Dinh}.

\begin{remark} \label{remarkDKW2} The significance of the set $\Eb$ appears most clearly in the ergodic theory of self-correspondences on curves over $\C$.  More precisely, for a balanced non-weakly modular self-correspondence, it is shown in \cite[Theorem 1.2]{DKW} that for every $x \in C(\C)$, $(\frac{1}{d_2}(\pi_1)_\ast (\pi_2)^\ast)^n(\delta_x) \rightarrow \mu_D$. Note that in the non-weakly modular case, it is easy to see that $\Kb=\Eb=\emptyset$.
When $d_1<d_2$, it is proved in \cite{Dinh} and \cite{DS} that $(\frac{1}{d_2}(\pi_1)_\ast (\pi_2)^\ast)^n(\delta_x) \rightarrow \mu_D$ for every $x \not \in \Eb$. It is therefore natural to conjecture that in every case $d_1 \leq d_2$, if $D$ is not finitary, $(\frac{1}{d_2}(\pi_1)_\ast (\pi_2)^\ast)^n(\delta_x) \rightarrow \mu_D$ if and only if $x \not \in \Eb$. At any rate, it is not hard to see that if $x \in \Eb$, if the limit $(\frac{1}{d_2}(\pi_1)_\ast (\pi_2)^\ast)^n(\delta_x)$ exists as a measure, then it charges at least one point in $\Kb$, and it cannot be $\mu_D$.
\end{remark}

\subsubsection{Polarized self-correspondences}

\begin{definition} Let $(D,\pi_1,\pi_2)$ be a self-correspondence on a curve $C$ over a field $k$. A {\it polarization} of $D$ is an ample line bundle $\Lc$ on $C$ such that $\pi_1^\ast \Lc^n = \pi_2^\ast \Lc^m$ for some positive integers $n$ and $m$. If $D$ admits a polarization we say that $D$ is {\it polarized}.
\end{definition}
Recall that on a curve, a line bundle is ample if and only if its degree is positive. If $\Lc$ is a polarization, one must have $d_1 n = d_2 m$.

A self-correspondence $D$ has a polarization in each of the following cases:
\begin{itemize} 
\item[(i)] $D=\P^1$ (in which case necessarily $C=\P^1$)
\item[(ii)] $k$ is algebraic over a finite field. 
\item[(iii)] $D$ has a finite non-empty equiramified complete set $S \subset C(k)$.
\end{itemize}
Indeed, let us consider the group homomorphism $h: \Pic\, C \rightarrow \Pic^0 D,\ \Lc \rightarrow (\pi_1^\ast \Lc)^{d_2} \otimes (\pi_2^\ast \Lc)^{-d_1}$. Clearly $D$ is polarizable if and only if $\ker h \not \subset \Pic^0 C$. In particular, $D$ is polarizable in case (i) since in this case $\Pic^0 D=0$ and $\ker h = \Pic C=\Z \not \subset \Pic^0 C=0$). Also $D$ is polarizable in case (ii) for in this case $\Pic^0 D$ is torsion while $\Pic\, C/\Pic^0 C =\Z$ is torsion free. In case (iii), if the self-correspondence $D$ has a finite equiramified  complete set $\emptyset \neq S \subset C(k)$, take $\Lc$ the line bundle attached to the divisor $\Delta_S = \sum_{s \in S} [s]$. 

\par \bigskip

We now recall the basics of the general height theory of Lang-Néron, following the exposition of  \cite[\S 3]{chambertloir} and \cite{serre}. Assume that $k$ is either a number field or a non-trivial finite type extension of an algebraically closed field. It is known that we can choose a family $M(k)$ of pairwise inequivalent absolute values on $k$ and numbers $\lambda_\nu > 0$ such that for $a \in k^\ast$, $|a|_\nu = 1$ for almost all $\nu \in M(k)$ and the product formula holds: $\prod_{\nu \in M(k)} |a|_\nu^{\lambda_\nu} =1$.
Moreover, if $k'$ is a finite extension of $k$, one can choose $M(k')$ in such a way that there  is a surjective map $sk') \rightarrow M(k)$ with finite fibers , such that for $a \in k$, $\nu \in M(k)$, $|a|_\nu^{\lambda_\nu} = \prod_{\nu' \in M(k'), \pi(\nu')=\nu} |a|_{\nu'}^{\lambda_{\nu'}}$, and the product formula hods.
The choice of $M(k)$ allows one to define the {\it height function} on $\P^n_{\bar k}$ by $$h([x_0,\dots,x_n]) = \log (\prod_{\nu \in M(k')} \max(|x_0|_\nu,\dots,|x_n|_\nu),$$ if $x_0,\dots,x_n \in k'$ (the result is independent of the choice of the finite extension $k'$ containing $x_0,\dots,x_n$) .  

For a $k$-projective variety $V$, let us denote by $\Fc(V)$ the $\R$-vector space of maps from $V(\bar k)$ to $\R$ and by $\Fcb(V)$ the subspace consisting of these maps that are bounded. There is a unique morphism  $\Pic\, V \rightarrow \Fc(V)/\Fcb(V)$, $\Lc \mapsto h_\Lc$, such that if $\Lc$ is very ample and $\phi$ is one of the embedding $V \rightarrow \P^n$ defined by $\Lc$, then $h_\Lc = h \circ \phi$. It follows that if $f:V \rightarrow W$ is a morphism of $k$-projective variety, $h_{f^\ast \Lc} =  h_\Lc \circ f$. 

\begin{lemma} If $\Lc$ is ample, then $h_{\Lc} \neq 0$ in $\Fc(V)/\Fcb(V)$. \end{lemma}
\begin{pf} We may assume that $\Lc$ is very ample, and this reduces us to prove that for a projective subvariety in $\P_k^n$, the function $h$ is unbounded on $C(\bar k)$. Up to a linear change of variables, the map $\pi([x_0,\dots,x_n])=[x_0,x_1]$ is surjective from $V$ to $\P^1$, and it is clear that $h(x) \geq h(\pi(x))$. It suffices therefore to show that $h$ is unbounded on $\P^1$, which is clear. 
\end{pf}
\begin{prop} \label{proppolarized} Let  $(D,\pi_1,\pi_2)$ be an unbalanced  polarized self-correspondence on a curve $C$ over a field $k$ that is not algebraic over a finite field. 
Then there are infinitely many vertices in $\Gamma_D$ that do not belong to any directed cycle.
\end{prop}
\begin{pf} 
Let $\Lc$ be a polarization on $D$. The self-correspondence $D$ together with $\Lc$  are defined over a subfield $k_0$ of $k$ which is of finite type over the prime subfield of $k$,
and if $k_1$ is any field such that $k_0 \subset k_1 \subset \bar k$ it suffices obviously to prove the result for $D$ considered as a self-correspondance over $k_1$. Since $\bar k$ is not the algebraic closure of a finite field, one can assume that $k_1$ is either 
a number field, or a non-trivial extension of finite type of an algebraic closed field. Replacing $k_1$ by $k$, we can now use the theory of height reminded above.

Let $(d_1,d_2)$
be the bi-degree of $D$. By symmetry of the statement to prove, we may and do assume that $d_1 < d_2$.

Since $\Lc$ is polarization, there exists an integer $n>0$ such that
$\pi_1^\ast(\Lc^{d_2 n}) = \pi_2^\ast(\Lc^{d_1 n})$. Thus $ d_2 h_{\Lc} \circ \pi_1= d_1 h_{ \Lc} \circ \pi_2$ in $\Fc(D)/\Fcb(D)$. In other word, if $h$ is any lift of $h_\Lc$ in $\Fc(C)$, there exists a positive real constant $M$ such that for all $z \in D(\bar k)$,
$|d_2 h(\pi_1(z)) - d_1 h (\pi_2(z))|<M$.

If $x_0,\dots,x_n$ is a directed cycle, then one has
\begin{eqnarray*} | h(x_0) - d_1/d_2 h(x_1)| &<& M/d_2 \\
\cdots & < & \cdots \\
| h(x_{n-1}) - d_1/d_2 h(x_0)| &<& M/d_2\end{eqnarray*}
so
$$|h(x_0) - (d_1/d_2)^n h(x_0) | <  \frac{M}{d_2 (1 + d_1/d_2 + \dots +(d_1/d_2)^{n-1})}< \frac{M}{d_2-d_1}.$$
It follows that if $x_0$ belongs to a directed cycle then $h(x_0)$ is bounded by a constant (independent of the length of the cycle).

By the lemma, $h$ is unbounded on $C(\bar k)$. There is therefore infinitely many points in $C(\bar k)$ that are not part of any directed cycle.
\end{pf}
\begin{remark} The proposition is obviously false for a finitary self-correspondence but we do not know whether the proposition holds for a polarized  balanced
non-finitary self-correspondence, nor whether the the polarized hypothesis may be dropped (even in the unbalanced case).
\end{remark}

\section{The operator attached to a self-correspondence}

\subsection{Definition of the operator $T_D$} \label{defTD}

If $C$ is a curve, and $(D,\pi_1,\pi_2)$ a self-correspondence of $C$, we denote by $T_{D,\pi_1,\pi_2}$ or simply $T_D$, the map $k(C) \rightarrow k(C)$ which sends $f$ to $$T_D f = \tr_{k(D) /_{\pi_2^\ast} k(C)} \pi_1^\ast(f).$$ 
The notations $\tr_{k(D) /_{\pi_2^\ast} k(C)}$ means the trace map from $k(D)$ to $k(C)$, where $k(D)$ is seen as an algebra over $k(C)$ of dimension $d_2$ through the map $\pi_2^\ast$.
The map $T_D$ is thus a $k$-linear endomorphism of $k(C)$.

\subsection{Local description of the operator $T_D$}

First recall some basic terminology. If $C/k$ is a curve, or a disjoint union of curves, 
and if $f \in k(C)$, $x \in C(\bar k)$, we denote by $\ord_x(f)$ the order of vanishing of $f$ at $x$. Thus $\ord_x(f) > 0$ if $f(x)=0$, $\ord_x(f)=0$ if $f(x) \in k^\ast$, and $\ord_x(f) < 0$ if $f(x) = \infty$. For $S$ a subset of $C(\bar k)$, we set $$\ord_S(f) := \inf_{x \in S} \ord_x(f).$$ If $S=C(\bar k)$, we simply write $\ord\, f$ for $\ord_S f$. One has $\ord\, f \leq 0$ for any $f$, and $\ord f = 0$ if and only if $f$ is a constant.

Now suppose given a self-correspondence $(D,\pi_1,\pi_2)$ of $C$ over $k$.

Given $y \in C(\bar k)$, and $z \in \pi_2^{-1}(y)$, we note $K_y$ and $K_z$ the fraction fields of the completions $\hat{\anneau_{C,y}}$ and $\hat{\anneau_{D,z}}$ of the local rings $\anneau_{C,y}$ and $\anneau_{D,z}$. The valuations $\ord_y$ on $k(C)$ (resp. $\ord_z$ on $k(D)$) extends uniquely to $K_y$ and $K_z$, and $\hat{\anneau_{C,y}}$ and $\hat{\anneau_{D,z}}$ are the rings of integers attached to these valuations; they are complete discrete valuation rings. The map $\pi_2^\ast$ induces an injective morphism $K_y \rightarrow K_z$, which makes $K_z$ a finite separable extension of $K_y$, totally ramified of degree $e_{2,z}$ (which means that 
$\ord_y(f) = e_2 \ord_z(f)$ for any  $f \in K_y$ seen as an element of $K_z$).
\begin{lemma} \label{lemmatrace}
Given $y \in C(\bar k)$, $z \in \pi_2^{-1}(y)$ and $f \in K_z$, let $P_{f}(X)=X^{d_2} + \sum_{i=1}^{d_2} a_{i,f} X^{d_2-i} \in K_y[X]$ be the characteristic polynomial of the multiplication by $f$ on $K_z$ seen as a $K_y$-vector space of dimension $d_2$ through $\pi_2^\ast$. Then one has, for $i=1,\dots,d_2-1$, $\ord_z(a_{i,f}) \geq i \ord_z(f)$ and $\ord_z(a_{d_2,f})={d_2} \ord_z(f)$. In particular,
$$\ord_y \tr_{K_z/K_y}(f) \geq  \lceil \ord_z(f) / e_{2,z} \rceil$$
\end{lemma} 
\begin{pf} 
Since $K_z/K_y$ is separable, $P_{f.z}(X)=\prod_{i=1}^{d_2} (X-\sigma_i(f))$ where the $\sigma_i$ runs amongst the embedding on $K_z$ into some normal closure $L$ of $K_Z$ over $K_y$. If $w$ is the valuation on $L$ extending $\ord_z$ on $K_Z$, then $w(\sigma_i(f))=\ord_z(f)$ and it follows that $\ord_z(a_{i,f})=w(a_{i,f}) \geq i w(f) = i \ord_z(f)$, with equality if $i=d_2$.
The last assertion follows since $a_{1,f}=\pm  \tr_{K_z/K_y}(f)$ and $\ord_z(a_{1,f}) = e_{2,z} \ord_y(a_{1,f})$.
\end{pf} 

\begin{prop} \label{local} Let $f \in k(C)$ and  $y \in C(\bar k)$. Set $$n = \min_{z \in D(k), \pi_2(z)=y}\left  \lceil \frac{e_{1,z} \ord_{\pi_1(z)} f} { e_{2,z}} \right \rceil.$$
Then $$\ord_{y} T_D f \geq n.$$ Moreover, if for any $z$ such that $\pi_2(z)=y$, $f$ has no pole at $\pi_1(z)$, then
$$f(y) = \sum_{z \in D(k), \pi_2(z)=y} e_{2,z} f(\pi_1(z)).$$
\end{prop}
\begin{pf} To compute the image of $T_D f$ in $K_y$ we may extend the scalars from $k(C)$ to $K_y$ since the formation of the trace commutes with base change. This means that $T_D f = \tr_{K_y \otimes_{k(C)} k(D)/ K_y} \pi_1^\ast(f)$ in $K_y$. But $K_y \otimes_{k(C)} k(D) = \prod_{z \in \pi_2^{-1}(y)} K_z$ (see e.g. \cite[Chapter II, \S3, Theorem 1]{serre2}) hence \begin{eqnarray} \label{TDfsum} T_D f = \sum_{z \in \pi_2^{-1}(y)} \tr_{K_z/K_y} \pi_1^\ast(f).\end{eqnarray}

For $z \in \pi_2^{-1}(y)$, setting $x=\pi_1(z)$, we have $\ord_z \pi_1^\ast f = e_{1,z} {\ord_x f}$ hence by Lemma~\ref{lemmatrace} $$\ord_y \tr_{K_z/K_y} \pi_1^\ast(f) \geq     \left  \lceil \frac{e_{1,z} \ord_{\pi_1(z)} f} { e_{2,z}} \right \rceil \geq n.$$ By (\ref{TDfsum}), $\ord_y T_D f \geq n$.

To prove the second assertion, note that under its assumption, for any $z$ such that $\pi_2(z)=y$, the image of $\pi_1^\ast f$ in $K_z$ belongs to the complete d.v.r. $\hat{\anneau_{D,z}}$ and its image in the residue field $\bar k$ is $f(\pi_1(z))$. Thus $\tr_{K_z/K_y} \pi_1^\ast(f) \in \hat{\anneau_{D,z}}$ and the image of that element in the residue field $\bar k$ is $e_{2,z} f(\pi_1(z))$. The formula $f(y) = \sum_{z \in D(k), \pi_2(z)=y} e_{2,z} f(\pi_1(z))$ then follows from (\ref{TDfsum})
\end{pf}

\begin{cor} \label{corTdf} Let $f \in k(C)$ seen as a map $f : C(\bar k) \rightarrow \bar k \cup \{\infty\}.$ Then the functions 
$T_D f$ and $A_{\Gamma_D} f$ agree on all points of $C(\bar k)$ but finitely many. They agree in particular on all étale complete sets 
on which $f$ has no pole. 
\end{cor}
\begin{pf} Indeed, the last formula of the above proposition shows that the two functions agree at all points $y \in C(\bar k)$ which are étale and not neighbors of a point where $f$ has a pole. 
\end{pf}

\begin{cor} if $D'$ and $D$ are two self-correspondence on $C$, $T_{D'D} = T_{D'} \circ T_D$.
\end{cor}
\begin{pf} For $f \in k(C)$, $T_{D'D} f$ agrees almost everywhere with $A_{\Gamma_{D'D}} f$, which agrees almost everywhere with  $A_{\Gamma_{D'}} A_{\Gamma_D} f$, which agrees almost everywhere with $T_D' T_D f$, and those two functions in $k(C)$ must then be equal.
\end{pf}

\subsection{The filtered ring $B_S$ attached to a set of vertices $S$}

Now  fix a self-correspondence $(D,\pi_1,\pi_2)$ on $C$ over $k$ and a finite non-empty complete set $S$ of $C(\bar k)$.
We denote by $B_S \subset k(C)$ the rings of rational functions on $C$ whose poles are all in $S$. Thus, $$B_S = \{f \in k(C),\ \ord_x(f) \geq 0 \text{ for all }x \in C(\bar k)-S\}.$$ if $S=\emptyset$, $B_S=k$. If $S$ is not empty, the ring of fractions of $B_S$ is $k(C)$.
For $n \geq 0$, we set $$B_{S,n} = \{f \in B_S,\ \ord_S(f) \geq -n\}.$$
\begin{lemma} The $k$-subspaces $B_{S,n}$ for $n = 0,1,\dots$ form an increasing exhaustive filtration of $B_S$. One has $B_{S,0}=k$. The quotients $B_{S,n}/B_{S,n-1}$ for $n \geq 1$ are spaces of dimensions $\leq |S|$, and of dimension exactly $|S|$ when $n$ is large enough.
\end{lemma}
\begin{pf}
The first two sentences are trivial. The last one follows from Riemann-Roch.
\end{pf}

\begin{prop} If $S$ is forward-complete, the subring $B_S$ of $k(C)$ is stable by $T_D$. If $S$ is forward-complete and ramification-increasing, the filtration $B_{S,n}$ is stable by $T_D$. The converses of both these statement hold if $\ch k = 0$ or $\ch k>d_2$.
\end{prop}
\begin{pf}
The first two statements follow from Prop.~\ref{local}. The converse statements are left to the reader.
\end{pf}

\begin{remark} Remember (Prop.~\ref{unbalancedprop}) that a forward-complete and ramification-increasing set $S$, non-empty and finite, may only exist if $d_1 \leq d_2$.  In the balanced case $d_1=d_2$, such a set $S$ has to be complete and equiramified (Scholium~\ref{scholiumbalanced}).
\end{remark}

\begin{example} If $C$ is the Igusa curve, and $S$ the super-singular complete set, which is étale, $B_S$ is the space of modular forms of level $N$, all weights, over $\bar \F_p$ and $(B_{S,n})$ is the {\it weight filtration} on that space  (see \cite{gross}). \end{example}

\subsection{Linearly finitary self-correspondences}

Let $D$ be a self-correspondence on a curve $C$ over $k$.
\begin{definition} We say that $D$ is {\it linearly finitary} if there is a non-zero polynomial $Q \in k[X]$ such that $Q(T_D)=0$ on $k(C)$.
\end{definition}

\begin{prop} \label{finimplin} If $D$ is finitary, then there is a monic polynomial $Q \in \Z[X]$ such that $Q(T_D)=0$ on $k(C)$. In particular $D$ is linearly finitary.
\end{prop}
\begin{pf} If $D$ is finitary, there exists an $M \geq 0$ such that every irreducible finite complete set of $D$ is a directed graph with $\leq M$ vertices, with $\leq d_1$ (resp. $\leq d_2$) arrows starting (ending) at each point. There are finitely many such graphs up to isomorphism, so infinitely many irreducible complete sets $S$ must be isomorphic to some finite directed graph $\Gamma$. If $Q(X)$ is the characteristic polynomial of the adjacency matrix of $\Gamma$, we see that for every $f \in k(C)$, $Q(T_D)f$ is zero on infinitely many irreducible complete sets, so $Q(T_D) f=0$, and therefore $Q(T_D)=0$.
\end{pf}

\begin{lemma} \label{RR}  Given two distinct points $p,q$ in $C(\bar k)$, a finite set $Z \subset C(\bar k)$ not containing $p$ or $q$, and an integer $n \geq 0$, there exists a rational function $f \in \bar k(C)$ such that $f(q)=1$, $f$ vanishes at every point of $Z$ at order at least $n$, and $f$ has no pole outside $p$.
\end{lemma}
This follows from Riemann-Roch.

\begin{lemma} \label{QTDzero} Let $Q(X) = \sum_{i=0}^n a_i X^i \in k[X]$ be a polynomial. The following are equivalent:
\begin{itemize}
\item[(i)] One has $Q(T_D )=0$ on $k(C)$.
\item[(ii)] There exists a non-empty forward-complete $S$ such that $Q(T_D)=0$ on $B_S$. 
\item[(iii)] There exists infinitely many irreducible étale complete sets $S'$ such that $Q(A_{\Gamma_D})=0$ on $\CC(S',k)$.
\item[(iii')] There exists infinitely many irreducible étale complete set $S'$ such that for every $x,x' \in S'$, one has $\sum_{i=0}^n a_i \np_{x,x',i}=0$ in $k$.
\item[(iv)] There exists an infinite étale complete set $S'$ such that $Q(A_{\Gamma_D})=0$ on $\CC(S',k)$.
\item[(iv')] There exists an infinite étale complete set $S'$ such that for every $x,x' \in S'$, one has $\sum_{i=0}^n a_i \np_{x,x',i}=0$ in $k$.
\item[(v)] There exists an infinite backward-complete set $S'$ such that $Q(A_{\Gamma_D})f$ has finite support for every  $f \in \CC(S',k)$.
\end{itemize}
\end{lemma}
\begin{pf} To prove that (i) implies (ii) is clear, take $S=C(\bar k)$.

To prove that (ii) implies (iii), let $p$ be a point in $S$. By Prop.~\ref{infinitecomplete}, there exists infinitely many étale complete sets $S'$ that do not contain $p$. Let $S'$ be one of then, and let $x$ be a point in $S'$. Let $\delta_x \in \CC(S',k)$ be the function whose value at $x$ is $1$ and is $0$ elsewhere.

We provide $\Gamma_D$ with the distance induced by its undirected graph structure.
By Lemma~\ref{RR}, and since $\Gamma_D$ is locally finite, there exists a function $f \in \bar k(C)$ such that $f=\delta_x$ on all points at distance $\leq 2n$ of $x$, and whose only possible pole is at $p$. In particular, $f \in B_S$, so $Q(T_D) f=0$. By Corollary~\ref{corTdf}, $Q(A_{\Gamma_D}) f =0 $ on $S'$. Clearly $Q(A_{\Gamma_D}) f$ and $Q(A_{\Gamma_D}) \delta_x$ agree on all points $x'$  at distance $\leq n$ of $x$. Since $Q(A_{\Gamma_D}) \delta_x$ has support in the sets of points at distance $\leq n$ of $x$, this implies that $Q(A_{\Gamma_D}) \delta_x = 0$. Since this is true for an arbitrary point $x$ of $S'$, $Q(A_{\Gamma_D})=0$ on $\CC(S',\bar k)$, hence (iii).

It is is clear that (iii) implies (iv), by taking $S'$ in (iv) to be the union al all the $S'$ in (iii). Also, the equivalences between (iii) and (iii') and (iv) and (iv') is just Formula~(\ref{Agn}).

Finally, (iv) implies (v) is clear, so it just remains to prove that (v) implies (i). Let $f \in k(C)$. Then $Q(A_{S'}) f$ and $T_D f$ agree on every points of $S'$ at distance $>n$ of a pole of $f$.
Since the number of such points is finite, $Q(A_{S'}) f$ and $T_D f$ agree on every points of $S'$ except a finite number of them, and 
this implies that $T_D f$ has finite support on $S'$, and since $S'$ is infinite, that $T_D f$ has infinitely many zeros. Hence $T_D f=0$.
\end{pf}

\begin{prop} If $Q(X) \in k[X]$ is a polynomial, then $Q(T_D)=0$ if and only if $Q(T_{{}^t D})=0$. In particular, $D$ is linearly finitary if and only if ${}^t D$ is. Moreover, if $k'$ is any field extension of $k$, then $D$ is linearly finitary if and only if $D_{k'}$ is linearly finitary.
\end{prop}
\begin{pf} All is clear using condition (iii') or (iv') of Lemma~\ref{QTDzero}.
\end{pf}

\begin{lemma} If $S$ is an infinite irreducible étale complete set in $\Gamma_D$, then for every integer $m$ and for every $x \in S$ there exists  $x' \in S$ with a directed path from $x$ to $x'$ of length $m+1$ and no directed path from $x$ to $x'$ of length $ \leq m$.
\end{lemma}
\begin{pf} By symmetry of the statement to be proved we may assume that $d_1 \geq d_2$ to begin with. For $x \in S$, denote by $F_x$ the smallest forward-complete set containing $x$, that is the set of all end points of directed paths starting at $x$. If $F_x$ is finite, then it is complete by Scholium~\ref{scholiumbalanced}, contradicting the irreducibility of $S$ which is infinite. Thus $F_x$ is infinite. 

Let $F_{x,m}$ be the subset of $F_x$ consisting of all end points of directed paths of length $\leq m$ starting at $x$. Then $F_{x,m} \subset F_{x,m+1}$ and $F_x = \cup_m F_{x,m}$. If for some $m$, $F_{x,m} = F_{x,m+1}$, then $F_{x,m}$ is forward-complete, hence $F_x=F_{x,m}$ and $F_x$ is finite, a contradiction. This for every $m$ there exists $x' \in F_{x,m+1}-F_{x,m}$, which proves the lemma.
\end{pf}

\begin{prop} \label{linfinfin} If $k$ has characteristic zero, and if $D$ is linearly finitary, then $D$ is finitary.
\end{prop}
\begin{pf} 
Assume $Q(T_D) = 0$ for $Q \in k[X]$ a non-zero polynomial of degree $n$ and dominant term $a_n \neq 0$. 
By Lemma~\ref{QTDzero}, there exists  infinitely many étale irreducible complete sets $S'$ such that $Q(A_{S'}) = 0$. For $S'$ any of them,  we prove by contradiction that $S'$ is finite. Indeed, choose $x \in S'$, and let $x'$ be a point in $S'$ with a directed path of length $n$ from $x$ to $x'$ but no shorter path, which exists by the lemma above if $S'$ is infinite; then $0=(Q(A_{S'}) \delta_x)=a_n \np_{n,x,x'}$, a contradiction since $\np_{n,x,x'} \geq 1$ is non-zero in $k$. Thus $D$ has infinitely many finite complete sets, and is therefore finitary.
 \end{pf}

\begin{cor} if $D$ is linearly finitary, there is a monic polynomial $Q \in \Z[X]$ such that $Q(T_D)=0$.
\end{cor}
\begin{pf} 
There is a non-zero monic polynomial $Q(X) = \sum a_i X^i \in k[X]$ such that $Q(T_D)=0$. Let $k_0$ be the prime subfield of $k$. Let $l$ be a $k_0$-linear form on $k$ such that $l(1) =1$ for some $i$. Then $Q_l(X)=\sum l(a_i) X^i \in k_0[X]$ is monic, and by the equivalence between (i) and (iii') in Lemma~\ref{QTDzero}, one has $Q_l(T_D)=0$. 
If $k_0$ is $\F_p$ for some $p$, it suffices to lift $Q_l$ into a monic polynomial in $\Z[X]$. If not, then $k$ has characteristic zero, so $D$ is finitary by Prop.~\ref{linfinfin}, and there exists a monic polynomial in $\Z[X]$ that kills $T_D$ by Prop~\ref{finimplin}.
\end{pf}

\subsection{The number of irreducible complete finite sets, II}

\label{finiteetaleII}

\subsubsection{Riemann-Roch calculations}

Let $k$ be an algebraic closed field, $C$ a curve of genus $g \geq 1$ (just to avoid modifying the formulas in the case $g=0$)
and $S \subset C(k)$ a finite set. We denote as above by $\Delta_S$ the effective Weil divisor $\sum_{s \in S} [s]$ in $\Pic\, C$.

Given a second finite non-empty set $S'$, disjoint from $S$, and any integer $n \geq 0$, there is clearly a unique relative integer $n'=n'(S,S',n)$ such that \begin{eqnarray} \label{defnprime}  2g-2 + |S'| \geq (n-1) |S|-n'|S'| >2g-2.\end{eqnarray} We denote by $V_{S,S',n}$ the subspace of $B_{S,n}$ of functions such that $\ord_{S'} f \geq n'$ where $n'$ is that integer, that is
\begin{eqnarray} \label{defV} V_{S,S',n}= \{ f \in k(C), \div f \geq - (n\Delta_S - n' \Delta_{S'}) \}.\end{eqnarray}

\begin{lemma} \label{BV} One has $B_{S,n-1}+V_{S,S',n} = B_{S,n}$ and $\dim V_{S,S',n}  \leq g-1 + |S|+ |S'|$.
\end{lemma}
\begin{pf}
The divisor $D =  (n-1) \Delta_S - n' \Delta_{S'}$ has degree $\deg D= (n-1) |S| - n' |S'| > 2g-2$ by assumption. Let $s_0 \in S$. By Riemann-Roch,  $h(D+[s_0])>h(D)$, so there exists a function $f_{s_0} \in k(C)$ such that $\div f_{s_0} \geq -D-[s_0]$ and $\ord_{s_0} f = -n$. Obviously $f_{s_0} \in V_{S,S',n}$ and the functions $f_{s_0}$, when $s_0$ runs in $S$, generate $B_{S,n}/B_{S,n-1}$. Hence the first assertion. 

The second assertion also follows from Riemann-Roch, which says that \begin{eqnarray*}
 \dim V_{S,S',n} &=& \deg (n \Delta_S - n' \Delta_{S'}) - g+1 \text{ (since $\deg (n \Delta_S - n' \Delta_{S'}) > 2g-2$ and $g>0$)} \\ &=& |S| + (n-1)|S| -n'|S'| -g+1\\ & \leq & |S| + 2g-2 + |S'| -g +1 \text{ (by (\ref{defnprime})} \\ &=& g-1 +|S| + |S'|. \end{eqnarray*}
\end{pf}

Now let $S''$ a third finite non-empty complete set, disjoint from $S$ and from $S'$. Let $n''$ be defined by (\ref{defnprime}) with $S'$ replaced by $S''$, and let $V_{S,S'',n'}$ be defined similarly with (\ref{defV}).
One has:
\begin{lemma} \label{interVSSn}There exists an integer $n_0$ such that for $n > n_0$, $V_{S,S',n} \cap V_{S,S'',n}=0$.
\end{lemma} 
\begin{pf} The space $V_{S,S',n} \cap V_{S,S'',n}$ is the space of functions $f$ such that $\div f \geq -(n\Delta S-n'\Delta_{S'}-n''\Delta_{S''})$.
Using (\ref{defnprime}) for $n'$ and $n''$, one computes:
\begin{eqnarray*} \deg( n\Delta_S-n'\Delta_{S'}-n''\Delta_{S''})&=&n|S|-n'|S'|-n''|S''|\\ &=&(n|S|-n'|S'|)+(n|S|-n''|S''|)-n|S|\\ & \leq & 2g-2 + |S'| + |S''| - n|S|.\end{eqnarray*} This number is negative for $n > (2g-2 + |S'| + |S''|)/|S|$, and therefore $V_{S,S',n} \cap V_{S,S'',n}=0$.
\end{pf}

\subsubsection{The bound}

\begin{theorem} \label{notthree} Let $D$ be a self-correspondence on a curve $C$ over an arbitrary field $k$. Assume that $T_D$ is not linearly finitary. Then $D$ has at most two irreducible complete equiramified finite sets.
\end{theorem}
\begin{pf}
We may and do assume that $k$ is algebraically closed. Also assume (just for simplicity in the formula) that the genus $g$ of $C$ is $\geq 1$, the case $C=\P^1$ being taken care of by Prop~\ref{caseP1}.
Assume that $D$ admits three irreducible complete equiramified finite sets, $S$, $S'$ and $S''$. By Prop.~\ref{local}, $V_{S,S',n}$
and $V_{S,S'',n}$ are stable by $T_D$. Let $n_0$ be as in Lemma~\ref{interVSSn}. We claim that for $n > n_0$, every eigenvalue $\lambda$ of $T_D$ in $B_{S,n}$ is also an eigenvalue of $T_D$ in $B_{S,n-1}$. We may assume that $\lambda$ is an eigenvalue of $T_D$ in the quotient $B_{S,n}/B_{S,n-1}$, otherwise there is nothing to prove. Let us call $m_\lambda \geq 1$ its multiplicity in $B_{S,n}/B_{S,n-1}$.  By Lemma~\ref{BV}, $\lambda$ is also an eigenvalue of $T_D$ in $V_{S,S'_n}$ (resp. in $V_{S,S'',n}$) with multiplicity $\geq m_\lambda$. Thus $\lambda$ appears as an eigenvalue of $T_D$ in $V_{S,S',n} + V_{S,S'',n}$ with multiplicity $\geq 2 m_\lambda$, since the sum is direct by Lemma~\ref{interVSSn}. Thus $\lambda$ appears as an eigenvalue of $T_D$ with multiplicity $\geq 2 m_\lambda >m_\lambda$ in $B_{S,n}$, and it must appear in $B_{S,n-1}$. 

By induction, all the eigenvalues of $T_D$ on $B_{S,n}$ (hence on $V_{S,S',n}$) for any $n$ already appear in $B_{S,n_0}$. Thus there are finitely many eigenvalues of $T_D$, say $\lambda_1,\dots, \lambda_l$, appearing in $V_{S,S',n}$ for any $n$.
Define $Q(X)=(X-\lambda_1)^{g-1+|S|+|S'|} \cdots (X-\lambda_l)^{g-1+|S|+|S'|}$. It is clear using Lemma~\ref{BV} that $Q(T_D)$ kills $V_{S,S',n}$ for any $n$, hence, by Lemma~\ref{BV}, $B_{S,n}$ for any $n$, and thus $Q(T_D)$ kills $B_S$, and by Lemma~\ref{QTDzero}, $Q(T_D)$ kills $k(C)$ contradicting the assumption that $D$ is not linearly finitary.
\end{pf}
\begin{cor} Let $D$ be a self-correspondence on a curve $C$ over a field $k$ of characteristic zero. Assume that $T_D$ is not finitary. Then $D$ has at most two irreducible complete equiramified finite sets.
\end{cor}
This follows from the theorem and Prop.~\ref{linfinfin}.

%
%
%

\end{document}